\documentclass[11pt]{amsart}
\usepackage[utf8]{inputenc}
\usepackage{amsmath,amssymb,proof,tikz}
\usetikzlibrary{matrix,arrows}
\newtheorem{theorem}{Theorem}[section] 
\newtheorem{prop}[theorem]{Proposition}
\newtheorem{lemma}[theorem]{Lemma}
\newtheorem{cor}[theorem]{Corollary}
\newtheorem{definition}[theorem]{Definition}
\newcommand{\rr}{\mathbb{R}}
\newcommand{\qq}{\mathbb{Q}}
\newcommand{\qp}{\mathbb{Q}_p}
\newcommand{\qqq}{\mathbb{Q}_q}
\newcommand{\zp}{\mathbb{Z}_p}
\newcommand{\cc}{\mathbb{C}}
\newcommand{\nn}{\mathbb{N}}
\newcommand{\zz}{\mathbb{Z}}

\newcommand{\nsr}{{^*\mathbb{R}}}

\newcommand{\nsn}{{^*\mathbb{N}}}

\newcommand{\al}{\mathcal{A}}
\newcommand{\almu}{\mathcal{A}_{\mu}}
\newcommand{\ral}{\mathcal{R}}
\newcommand{\sal}{\mathcal{S}}
\newcommand{\ralmu}{\mathcal{R}_{\mu}}

\newcommand{\st}{\text{\textit{st}}\hspace{0.2mm} }

\newcommand{\ns}{ \text{\textit{ns}}\hspace{0.5mm} }
\newcommand{\pns}{\text{\textit{pns}}\hspace{0.5mm} }
\newcommand{\fin}{\text{\textit{fin}}\hspace{0.5mm} }
\newcommand{\infi}{\text{\textit{inf}}\hspace{0.5mm} }
\newcommand{\mon}{\text{\textit{mon}}\hspace{0.5mm} }

\begin{document}
 
\author{Heiko Knospe }
\address{Technische Hochschule K\"oln, Faculty 07, Betzdorfer Str. 2, D 50679 Köln, Germany} 
\email{heiko.knospe@th-koeln.de}
\title{Nonstandard Measure Spaces with Values in non-Archimedean Fields}

\begin{abstract}
The aim of this contribution is to bring together the areas of $p$-adic analysis and nonstandard analysis. We develop a nonstandard measure theory with values in a complete non-Archimedean valued field $K$, e.g.\ the $p-$adic numbers $\mathbb{Q}_p$. The corresponding theory for real-valued measures is well known by the work of P. A. Loeb, R. M. Anderson and others. 

We first review some of the standard facts on non-Archimedean measures and briefly sketch the prerequisites from nonstandard analysis. Then internal measures on rings and algebras with values in a nonstandard field ${^*K}$ are introduced. We explain how an internal measure induces a  $K$-valued Loeb measure. 
The standard-part map between a Loeb space and the underlying standard measure space is measurable almost everywhere. We establish liftings from measurable functions to internal simple functions. Furthermore, we prove that standard measure spaces can be described as push-downs of hyperfinite internal measure spaces. This result is an analogue of a well-known Theorem on hyperfinite representations of Radon spaces. Then standard integrable functions are related to internal $S$-integrable functions and integrals are represented by hyperfinite sums. Finally, the results are applied to measures and integrals on $\mathbb{Z}_p$ and $\mathbb{Z}_p^{\times}$. We obtain explicit series expansions for the $p$-adic zeta function and the $p$-adic Euler-Mascheroni constant which we use for computations.

\end{abstract}

\maketitle
\tableofcontents




\section{Introduction}

Integrals of functions with values in a complete non-Archimedean field are studied in the field of {\em $p$-adic analysis} and a general measure-theoretical approach to $p$-adic integration has been developed by A. van Rooij \cite{rooij}. 
$p$-adic measures and integrals are used in number theory and in arithmetic geometry, in particular in the context of $p$-adic zeta- and $L$-functions. 

This contribution applies methods of {\em nonstandard analysis} to measures and functions with values in a complete non-Archimedean field, e.g.\ the $p$-adic numbers $\qq_p$. {Nonstandard \mbox{analysis}} was established   in the 1960s by A. Robinson \cite{robinson1961},\cite{robinson1996}.  
Nonstandard extensions can be defined by the {\em ultrapower} construction and they behave in a functorial way \cite{serpe}. 
In the past, nonstandard analysis has been successfully applied to real measure theory. {\em Loeb measures} \cite{albeverio1986} are of particular importance because they permit the transition from nonstandard to standard measure spaces. Our aim is to investigate $p$-adic measure spaces with nonstandard methods and to obtain Theorems similar to the case of real Radon measures, for example representations by hyperfinite measure spaces. \\

\subsection{Measures with values in non-Archimedean fields}
\label{measures}
Let $K$ be a field with a non-Archimedean absolute value $|\ |$. We suppose that the absolute value is non-trivial and $K$ is complete. 
We recall some basic definitions and facts on $K$-valued measures and integrals from \cite{rooij}, where a detailed exposition of the subject can be found. Our description is based on measures and measurable functions rather than distributions and continuous functions. Of course, these terms and definitions are closely related.

Let $X$ be a set and $\ral$ a {\em ring} of subsets of $X$. This means that $\varnothing \in \ral$ and for any $A, B \in \ral$, we have $A \cup B \in \ral$ and $A \setminus B \in \ral$. We assume $\ral$ is {\em covering} and {\em separating}, i.e., for any $a, b \in X$ there is $A \in \ral$ such that $a\in A$ and $b\in X \setminus A$. The sets in $\ral$ are called {\em measurable} and $(X,\ral)$ is called a {\em measurable space}. $\ral$ is the base of a {\em zero-dimensional} Hausdorff topology on $X$. If in addition $X \in \ral$,  then $\ral$ is an {\em algebra} and we will see that measures on algebras have additional favourable properties.\\
 
For a given zero-dimensional Hausdorff space $X$, let $B(X)$ be the set of clopen (open and closed) subsets. Then $B(X)$ is a covering and separating algebra.
 For a locally compact space $X$, the compact and clopen sets form a covering and separating subring of $B(X)$ which is denoted by $B_c(X)$. For example, if $K$ is a locally compact field, then $B_c(K)$ consists of all finite unions of bounded balls $B_r(a) = \{ x \in K\ |\ |x-a| \leq r \}$.

\begin{definition} Let $(X,\ral)$ be a measurable space. A {measure}  on $\ral$ with values in $K$ is a map $\mu: \ral \rightarrow K$ with the following properties:
\begin{enumerate}
\item $\mu(A \cup B)=\mu(A)+\mu(B)$ for disjoint sets $A, B \in \ral$. (Additivity)
\item For $A \in \ral$, $\|A\|_{\mu} = \sup\ \{|\mu(B)|\ :\ B \subset A,\ B\in \ral \} < \infty$. \\
(Boundedness)
\item For any shrinking set $\ral_0 \subset \ral$ (i.e., $A, B \in \ral_0 \Rightarrow A \cap B \in \ral_0$) with empty intersection (i.e., $\bigcap_{A \in \ral_0} A = \varnothing$) and for any $\epsilon > 0$ there exists a set $A_{\epsilon} \in \ral_0$ with $|\mu(A_{\epsilon})|<\epsilon$. (Continuity)
\end{enumerate}
$(X,\ral,\mu)$ is called a measure space.
\label{measure}
\end{definition}

\noindent {\em Remarks.} If one of the sets $A \in \ral_0$ in part (c) is compact in the $\ral$-topology, then c) is automatically satisfied, since in this case
a finite intersection must be empty.
Furthermore, the continuity implies $\sigma$-additivity, i.e.\ $\mu(\bigcup_{n\in \nn} A_n) = \sum_{n\in \nn} \mu(A_n)$ for disjoint sets $A_n \in \ral$, if $A=\bigcup_{n\in \nn} A_n \in \ral$ holds. In fact, the sequence $A$, $A \setminus A_1$, $A \setminus (A_1 \cup A_2)$, $\dots$ is shrinking and therefore c) yields $\lim_{n \rightarrow \infty} \mu(A)-(\mu(A_1)+\dots+\mu(A_n)) = 0$.
But note that is not required that infinite unions are measurable and  $\ral$ is usually not a $\sigma$-algebra (even if $\ral$ is an algebra).
The reason for this difference to real measure spaces is that a $K$-valued measure on a $\sigma$-algebra is purely atomic and therefore almost trivial (see \cite{rooij} 4.19 and 7.A): suppose for example that $X$ is a compact ultrametric with the $\sigma$-algebra of Borel sets; then the singletons $\{a\}$ are measurable and the continuity property implies that the measure values of a shrinking set of punctured balls with a fixed center converges to $0$. Any measurable set $Y$ can be covered by a finite set of balls where the absolute value of the measure of the punctured balls is small, i.e.\ less than any fixed $\epsilon$. Since the absolute value of $K$ is non-Archimedean, one obtains $|\mu(Y)| \leq \epsilon$ and hence $\mu(Y)=0$, unless the measure is atomic.

The real number $\|A\|_{\mu}$ is
defined as the supremum of all $|\mu(B)|$ where $B$ is a subset of $A$.
If $\ral$ is an algebra then $\|X\|_{\mu}$ is a global upper bound for all $|\mu(A)|$. 
The definition implies $\|A\|_{\mu} \leq \|B\|_{\mu}$ for $A, B \in \ral$ and $A \subset B$, as one expects. 
The reason for considering $\|.\|_{\mu}$ in addition to $|\mu(.)|$ is that the latter is not monotone; $| \mu(A) | \leq |\mu(B)|$ can be false for a subset $A \subset B$. But the inequality $|\mu(B)| \leq \max(|\mu(A)|, |\mu(B\setminus A)|)$ holds. If $A,B \in \ral$ then $\|A \cup B\|_{\mu} \leq \max \{ \|A\|_{\mu}, \|B\|_{\mu} \}$. Indeed, if $C \subset A \cup B$, then $\mu(C) = \mu(C \cap A) + \mu(C \cap (B \setminus A))$ and hence 
$$|\mu(C)|_{\mu} \leq \|C \cap A\|_{\mu} + \|C \cap (B \setminus A)\|_{\mu} \leq \|A\|_{\mu} + \|B\|_{\mu}. $$

One can show that part (c) of the above Definition implies the existence of a set $A_{\epsilon} \in \ral_0$ with the stronger property $\|A_{\epsilon}\|_{\mu}<\epsilon$ (see \cite{rooij} chapter 7). Thus (c) is equivalent to the following:\\

\begin{enumerate}
\item[(c')] For any shrinking set $\ral_0 \subset \ral$  with empty intersection one has $\displaystyle\lim_{A \in \ral_0} \| A\|_{\mu} = 0$.
\end{enumerate}

We say $A \in \ral$ is a {$\mu$-\em null set} or {$\mu$-\em negligible} if $\|A\|_{\mu} = 0$. This notion will be extended below to arbitrary subsets of $X$.
 One defines a {\em Norm function} $N_{\mu} : X \rightarrow \rr_{\geq 0}$ by
$$N_{\mu}(x) = \inf \{\|A\|_{\mu}  : x\in A \in \ral \}$$
 $N_{\mu}$ is called a Norm function, since it is used to define a seminorm  on the space of $K$-valued functions on $X$ (see Section \ref{sec:int}).
 For real-valued regular measures, such a function $N_{\mu}$ would mostly be zero (except for atomic measures), but this is not the case for $K$-valued measures. In fact, for every locally compact space $X$ there exists a measure $\mu$ on $B_c(X)$ such that $N_{\mu}=1$ everywhere (\cite{rooij} 7.9).\\
 
$\|A\|_{\mu}$  can be recovered from $N_{\mu}$ by the formula  $\|A\|_{\mu} = \sup_{x\in A} N_{\mu}(x)$ (see \cite{rooij} 7.2). The inequality $\geq$ follows from the definition of $N_{\mu}$. For the reverse inequality, take $\epsilon >0$. For every $x \in A$ there exists $B \in \ral$ with $x \in B$ such that $\|B\|_{\mu} \leq N_{\mu}(x)+\epsilon$. Using the continuity of the measure (the complements $A \setminus B$ of such $B$'s are shrinking with empty intersection), one finds a $B \in \ral$ with  $\|B\|_{\mu} \leq (\sup_{x \in A} N_{\mu}(x))+\epsilon$ and $\|A \setminus B\|_{\mu} \leq \epsilon$. Thus $\|A\|_{\mu} \leq \max\{\|B\|_{\mu},\| A\setminus B\|_{\mu} \} \leq
 (\sup_{x\in A} N_{\mu}(x))+\epsilon$ which establishes the asserted inequality.\\

It follows that $A \in \ral$ is a $\mu$-null set if and only if $N_{\mu}(x)=0$ for all $x \in A$. The latter is used to define $\mu$-null subsets of $A$ which are not necessarily measurable. Since null sets are thus compatible with arbitrary unions, there exists a largest $\mu$-null subset of $X$. \\

Now we extend our ring $\ral$ and include all sets which can be approximated by measurable sets. 
 The extended (completed) ring $\ralmu \supset \ral$ contains all $A \subset X$ with the following property: for all $\epsilon > 0$ there is a measurable set $B_{\epsilon} \in \ral$ such that $N_{\mu}(x) \leq \epsilon$ for all $x \in A \Delta B_{\epsilon} = (A \cup B_{\epsilon}) \setminus (A \cap B_{\epsilon})$. The latter is the symmetric difference of $A$ and $B_{\epsilon}$. This means that a set $A \in \ralmu$ can be approximated by a set $B_{\epsilon} \in \ral$. If $B_{\epsilon}' \in \ral$ is another approximating set then $N_{\mu}(x) \leq \epsilon$ for all $x \in B_{\epsilon} \Delta B_{\epsilon}'$. This implies $|\mu(B_{\epsilon})-\mu(B'_{\epsilon})| \leq \|B_{\epsilon} \Delta B_{\epsilon}'\|_{\mu} \leq \epsilon $. 
It can be easily shown that $\ralmu$ is again a ring. The measure $\mu$ can be extended to $\ralmu$ by taking the limit $ \mu(A)= \lim_{\epsilon \rightarrow 0} \mu(B_{\epsilon})$ which is well defined by the above. The additivity and boundedness is obvious. Let $X_{\epsilon} = \{ x\in X : N_{\mu}(x)>\epsilon \}$. The continuity of the extended measure follows by intersecting the elements of a shrinking subset with $X_{\epsilon}$ and using the continuity of the original measure (see \cite{rooij} 7.4).
The extended algebra is {\em complete}, i.e.\ it contains all subsets of $\mu$-null sets. 
We obtain the {\em extended measure space} $(X,\ralmu,\mu)$. $\ralmu$ is the base of the zero-dimensional $\ralmu$-topology on $X$ which is finer than the original $\ral$-topology. Furthermore, $\ralmu$ is stable under further extension, i.e., a set which can be approximated by sets in $\ralmu$ is already contained in $\ralmu$ (see \cite{rooij} 7.5). \\

The following statement (\cite{rooij} 7.6) is crucial for the next section:

\begin{prop} Let  $(X,\ralmu,\mu)$ an extended measure space, $A \in \ralmu$, $\epsilon>0$ and $X_{\epsilon} = \{ x\in X : N_{\mu}(x)>\epsilon \}$. Then $X_{\epsilon} \cap A$ is $\ralmu$-compact and $N_{\mu}$ is $\ralmu$-upper semicontinuous. 
\label{compact}
\end{prop}

We obtain an interesting relationship between measurable sets and its associated \mbox{topology}:

\begin{cor} Let $B$ be any clopen set in $X$ with respect to the zero-dimensional $\ralmu$-topology and suppose that $B \subset A$ for some $A \in \ralmu$. Then $B \in \ralmu$. In particular, if $\ral$ is an algebra and $B(X)$ the algebra of clopen sets in the $\ralmu$-topology, then $B(X)=\ralmu$.
\label{bx}
\end{cor}
{\em Proof.} Let $B \in B(X)$ and $\epsilon>0$. Then Proposition \ref{compact} implies that $A \cap X_{\epsilon}$ and hence also $B \cap X_{\epsilon}$ is compact. Since $B$ can be covered by sets in $\ralmu$, there is a finite sub-cover of $B \cap X_{\epsilon}$ which yields $B \cap X_{\epsilon} \in \ralmu$. Since this holds for all  $\epsilon>0$, we then have $B \in \ralmu$. $\hfill \square$\\

{\noindent \em Remark.} A measure $\mu$ on the algebra $B(X)=\ralmu$ is called {\em tight}. In this case,  $X_{\epsilon}$ is compact and $\| X \setminus X_{\epsilon}\|_{\mu} \leq \epsilon$ for all $\epsilon > 0$. \\

There is a close connection between measurability and continuity:

\begin{cor} The following conditions are equivalent:
\begin{enumerate}
\item $f: X \rightarrow K$ is $(\ralmu, B(K))$-locally measurable, i.e. $f \cdot \chi_A$ is measurable for any $A \in \ralmu$.
\item $f: X \rightarrow K$ is $(\ralmu, B(K))$-continuous. 
\end{enumerate}
\label{cont}

\end{cor}
{\em Proof.} a) implies b) since $\ralmu$ is a covering ring and conversely, b) implies a) by the above Corollary \ref{bx}.$\hfill\square$ \\

{\noindent \em Remark.} If $\ral$ is only a ring, then even the constant functions are only locally measurable. If $\ral$ is an algebra then {\em locally measurable} can be replaced by {\em measurable} in part a).\\

If $N_{\mu}$ is everywhere greater than some positive number, then $\ral=\ralmu$:

\begin{cor} Let $(X,\ral,\mu)$ be a measure space as above and assume that there is some $\epsilon>0$ such that $N_{\mu}(x) > \epsilon$ for all $x\in X$. Then $\ral=\ralmu$, $X$ is locally compact in the $\ral$-topology and
a function $f: X \rightarrow K$ is $(\ral, B(K))$-continuous if and only if $f$ is $(\ral, B(K))$-locally measurable. 
\end{cor}

\subsection{Nonstandard extensions}
\label{nonst}
In this subsection, we briefly recall the notions and prerequisites from nonstandard analysis (see for instance \cite{albeverio1986}, \cite{lr}, \cite{vaeth} for more details). In particular, we discuss nonstandard interpretations of $p$-adic fields. 

Nonstandard analysis was invented by 
 by Abraham Robinson in 1961 \cite{robinson1961}. His original construction uses model theory and was motivated by a Theorem of T. Skolem who showed the existence of {\em nonstandard} models of arithmetic: the natural numbers can not be uniquely characterized with first-order logic \cite{skolem1934}. There are countable models with infinite (unlimited) numbers. Similarly, there exists a nonstandard model $\nsr$ of the theory of the real numbers. $\nsr$ is an ordered extension field of $\rr$ which contains numbers greater than any standard real number \cite{robinson1996}. Since $\nsr$ is a field, it also contains infinitesimal numbers. Later W.A.J. Luxemburg gave an explicit construction of the hyperreal numbers by equivalence classes of sequences of real numbers modulo an ultrafilter (the {\em ultrapower} construction) which is widely used today.\\ 
 
This construction can be applied to almost any mathematical object which is contained in a {\em superstructure} $V(S)$ above some base set $S$. The latter is obtained by iterating the power set operation over the base set (for example $S=\rr$) and taking unions. There is a  general embedding map $^* : V(\rr) \rightarrow V(\nsr)$ called {\em nonstandard extension} between the superstructures over $\rr$ and $\nsr$. An object $A \in V(\rr)$ (e.g., a set like $\nn$ but also higher order structures such as fields, functions, topological spaces or measure spaces) is mapped to an extended object ${^* A} \in V(\nsr)$. ${^*A}$ can be defined by the ultrapower construction, similar to the case of the real numbers.  Taking ultrapowers over the index set $I=\nn$ yields a  countable saturated embedding. Countable saturation suffices for many applications (e.g.\ for countable sets and separable spaces), but sometimes richer nonstandard embeddings are required, i.e.\ a $\kappa$-saturated or polysaturated embedding of superstructures (see \cite{vaeth} for more details). 
The elements $A \in V(\rr)$ are called {\em standard}.
Let $B \in V(\nsr)$. Then $B$ is called {\em standard} if $B={^*A}$, i.e.\ $B$ is obtained by a constant sequence of $A$'s. $B$ is called {\em internal} if $B \in {^*A}$, i.e.\ $B$ is represented by a sequence $(a_i)_{i\in I}$ with $a_i \in A$. The remaining objects in $V(\nsr)$ are called {\em external}. Note that the standard copy of $A \in V(\rr)$ in $V(\nsr)$, i.e.\ the set $^{\sigma}A=\{ {^*a} \in {^*A} \ |\ a \in A\}$, is an external subset of the standard set ${^*A}$. 

\noindent {\em Example:} $\nsn$ can be constructed as the product of copies of $\nn$ modulo the given ultrafilter. 
$\nsn$ contains infinite numbers, for example the class of the sequence $(n)_{n\in\nn}$.
A sequence of finite subsets of $\nn$, i.e.\ a sequence of elements in $A=\mathcal{P}(\nn)$, gives an element of ${^*\mathcal{P}(\nn)}$ and hence an internal subset of $\nsn$. This internal set is standard-finite if it coincides with some fixed finite set on a set of indices contained in the ultrafilter. Otherwise, the set is hyperfinite, i.e.\ internal and nonstandard. If $N \in \nsn$ is an infinite number then $\{1,2,\dots, N\}$ is a hyperfinite set which contains a copy of $\nn$. Note that the subsets $^{\sigma}\nn$ and  $\nsn \setminus {^\sigma}\nn$ of $\nsn$ are external.\\

The embedding ${^*}:\ \mathcal{C} \rightarrow {^* \mathcal {C}}$ of a small category $\mathcal{C}$ (relative to the superstructure) is a covariant functor, and a given functor $F  :\mathcal{C}_1 \rightarrow \mathcal{C}_2 $ between small categories can be extended to a functor ${^* F} : {^* \mathcal {C}_1} \rightarrow {^* \mathcal {C}_2}$. It can be shown that these functors are well behaved (see \cite{serpe}).
There are a number of important principles which we mention only shortly:
\begin{enumerate}
\item {\em Transfer Principle}: Terms, formulas and sentences can be extended to the nonstandard universum. Objects are extended by the $*$-embedding and quantifiers over sets are substituted by quantifiers over internal sets. The important transfer principle states that a sentence $\varphi$ is true if and only if $^* \varphi$ is true.
\item {\em Saturation Principle}: Suppose that a family of internal sets $(A_i)_{i \in I}$ has nonempty finite intersections, the nonstandard embedding $^*$ is $\kappa$-saturated and the cardinality of $I$ is at most $\kappa$. Then $\bigcap_{i \in I} A_i \neq \varnothing$.
\item  {\em Countable Saturation}: A countable decreasing sequence $(A_n)_{n\in \nn}$ of  nonempty internal sets has a nonempty intersection $\bigcap_{n\in \nn} A_n$.
\item {\em Countable Comprehension} (equivalent to countable saturation): A sequence $(a_n)_{n\in \nn}$ of elements of an internal set $A$ can be extended to an internal sequence $(a_n)_{n\in \nsn}$ of elements of $A$. 
\item {\em Permanence Principle}: Let $A(n)$ an internal formula with $n$ the only free variable. If $A(n)$ holds for all $n \in \nn$ with $n \geq n_0$, then there exists an infinite $N_0 \in \nsn$ such that $A(n)$ holds for all $n \in \nsn$  with $n_0 \leq n \leq N_0$. If $A(n)$ holds for all infinite $n \in \nsn$, then there exists $n_0 \in \nn$ such that $A(n)$ holds for all $n \geq n_0$.
\end{enumerate} 

A topological space $(X,\mathcal{T})$ possesses an extension $({^* X},{^*\mathcal{T}})$. Let $a \in X$. Then the intersection of all standard  neighbourhoods  of $a$ in ${^* X}$ is called the {\em monad} of $a$: $$\mon(a) = \bigcap_{a \in A \in \mathcal{T}} {^*A}$$
 An element in $x \in {^*X}$ is called {\em nearstandard} if $x \in mon(a)$ for some $a\in X$ and $\ns({^*X})$ is the subset of nearstandard elements. 
 If $X$ is a Hausdorff space, then $mon(a) \cap mon(b)=\varnothing$ for $a \neq b$. This allows to define the important {\em standard-part} function $\st_X : \ns({^*X}) \rightarrow X$ which projects the subset $mon(a)$ to $a$. To simplify notation we often write $\st$ instead of $st_X$. 

The following statement (see \cite{lr} 21.7) gives a nonstandard characterisation of open, closed and compact sets. We remark that countable saturation of the superstructure embedding suffices if the topological space is separable.

\begin{prop} Let $X$ be a Hausdorff space and $A \subset X$. Then:
\begin{enumerate}
\item $A$ is open if and only if $\st^{-1}(A) \subset {^* A}$, i.e.\ all monads are contained in ${^*A}$.
\item $A$ is closed if and only if ${^*A}\, \cap\, \ns({^*X}) \subset \st^{-1}(A)$, i.e.\ the nearstandard elements in ${^*A}$ are contained in some monad of $A$.
\item $A$ is compact if and only if ${^*A} \subset \st^{-1}(A)$, i.e.\ all elements in ${^*A}$ are contained in some monad of $A$. This is equivalent to $\ns({^*A})={^* A}$.
\end{enumerate}
\label{topol}
\end{prop}

Since  $\st^{-1}(A) \subset \ns({^* A})$ holds by definition, we obtain the following Corollary:

\begin{cor} Let $X$ be a Hausdorff space. Then $A \subset X$ is clopen if and only if $\st^{-1}(A) = {^* A} \cap \ns({^*X})$.
\label{clopen}
\end{cor}   

We next turn to nonstandard extensions of fields. Let $K$ be a field which is complete with respect to a non-Archimedean absolute value $|\ |_v$. Let $o_K=\{x \in K\ :\ |x|_v \leq 1 \}$ be the ring of integers, $\frak{m}_K=\{x \in K\ :\ |x|_v <1 \}$ its maximal ideal and $k=o_K/\frak{m}_K$ the residue field of $K$. 
The absolute value $|\ |_v$ extends to an internal absolute value $|\ |_{^*v}: {^*K} \rightarrow \nsr_{\geq 0}$. We write $|\ |$ for either $|\ |_v$ or $|\ |_{^*v}$ for simplicity of notation. 
An element $x \in {^* K}$ is called {\em finite}, if $|x|$ is a finite hyperreal number and is called {\em infinitesimal} if $\st_K(x)=0$ or equivalently  $\st_{\rr} |x| = 0$. The set $\fin({^*K})$ of finite elements is a {\em valuation ring} which includes $K$. The ideal $\infi({^*K})$ of infinitesimal elements is the maximal ideal of the valuation ring $\fin({^*K})$ and their residue field is studied below. 
Two elements $x,y \in {^*K}$ are called {\em approximate} ($x \approx y$) if their difference is infinitesimal. The set $x+\infi({^*K})$ of elements which are approximate to $x \in K$ is the {monad} of $x$ and the  union of all monads is the subset $\ns({^*K})$ of nearstandard elements.

The internal absolute value defines a uniform structure and a Hausdorff topology on ${^*K}$. 
The topology is the nonstandard extension of the topology on $K$ and ${^*K}$ then has the structure of a topological field.
It is easy to see that $\infi({^*K})$, all monads,  $\ns({^*K})$ and $\fin({^*K})$ are clopen subsets. The concatenation $\st_{\rr} \circ |\ |_v$ defines a {\em seminorm} on $\fin({^*K})$; the corresponding topology is not Hausdorff and coarser than the topology defined by the internal absolute value.
The seminorm on $\fin({^*K})$ induces a norm on the quotient space $ \fin({^*K}) / \infi({^*K})$.

An element $y \in {^* K}$ is called {\em pre-nearstandard} ($y \in \pns({^*K})$) if for any standard $\epsilon>0$ there is some $x \in K$ such that $|x-y| < \epsilon$ (cf. \cite{lr} 24.7).
There are obvious inclusions
\[ \infi({^*K}) \subset \ns({^*K})  \subset \pns({^*K}) \subset \fin({^*K}) \subset {^*K} . \]

\begin{prop} Let $K$ be a locally compact non-Archimedean field. Then
$\ns({^*K}) = \fin({^*K})$ .
\label{pnsfin}
\end{prop}
{\em Proof.} Let $x \in {^*K}$ be finite. Then $x$ is contained in some open unit ball ${^*B}_r(0)$ with standard radius $r$. For any $\epsilon>0$, $B_r(0)$ can be covered with balls $B_{\epsilon}(y)$ with radius $\epsilon>0$ and some center $y \in K$. Since $B_r(0)$ is compact, $B_r(0)$ can be covered with a {\em finite} number of balls $B_{\epsilon}(y)$. Since the $^*$-embedding commutes with finite unions, we obtain $x \in {^*B}_{\epsilon}(y)$ so that $x$ can be approximated by $y \in K$. Hence $x$ is pre-nearstandard. Furthermore, the completeness of $K$ implies $\ns({^*K})=\pns({^*K})$ (see \cite{lr} 24.15) which gives our assertion.$\hfill\square$\\

\begin{prop} The standard-part map induces an isometric isomorphism $ \ns({^*K}) / \infi({^*K}) \cong K $. 
For a locally compact field $K$ one has \\ $ \fin({^*K}) / \infi({^*K}) \cong K$.
\label{Kv}
\end{prop}
{\em Proof.} The standard-part maps commutes with the norm: $st_{\rr} ( |x|_v) = |\st_K(x)|_v$ for $x \in \fin({^*K})$. The isomorphism follows from a general Theorem on nonstandard representations of complete metric spaces (see \cite{lr} 24.19). For locally compact fields one applies Proposition \ref{pnsfin}.$\hfill\square$\\

\noindent {\em Remark.} The above Proposition also holds for non-complete fields $K$, if one replaces $ \ns({^*K})$ by $ \pns({^*K})$ and completes $K$ on the right-hand side. This construction  can be used to define $\cc_p$   and a spherical completion $\Omega_p$. Since $\overline{\qp}$ is not locally compact, the $\pns$ subspace of ${^*\overline{\qp}}$ is a strict subset of the the $\fin$ subspace.  \\

The following Lemma is easy to show:

\begin{lemma} Let $(a_n)_{n\in \nn}$ be a sequence of elements in $K$. It extends to a sequence $(a_n)_{n\in \nsn}$ in ${^*K}$ and
\begin{enumerate}
\item $\lim_{n\rightarrow \infty} a_n = a \in K$ if and only if $a_N \approx a$ for all infinite $N \in \nsn$.
\item The series $\sum_{n=0}^{\infty} a_n $ converges in $K$ if and only if  $a_N \approx 0$, or equivalently $|a_N| \approx 0$, for all infinite $N \in \nsn$. 

\end{enumerate}
\label{convergence}
\end{lemma}

\section{Measure spaces}

\subsection{Internal measure spaces}
In this section, we define internal measure spaces with values in non-Archimedean fields.  Let $\Omega$ be an internal set for a given nonstandard extension (e.g. $\Omega={^*X}$ for a standard set $X$) and $\sal$ an internal covering and separating subring of $\mathcal{P}(\Omega)$ (e.g. 
$\sal = {^*\ral}$ for a standard ring $\ral$). Then $(\Omega,\sal)$ is a measurable space. Let $K$ be a complete non-Archimedean valued field and $B(K)$ the rings of clopen subsets of $K$. Then $(K,B(K))$ and $({^*K},{^*B(K)})$ are measurable spaces.

\begin{definition} An internal measure on $\sal$ with values in $^*K$ is an internal function $\nu: \sal \rightarrow {^*K}$ such that
\begin{enumerate}
\item $\nu(A \cup B)=\nu(A)+\nu(B)$ for disjoint sets $A, B \in \sal$. (Additivity)
\item For all $A \in \sal$, $\|A\|_{\nu} = \sup\ \{|\nu(B)|\ :\ B \subset A,\ B\in \sal \} $ is a hyperreal number. (Boundedness)
\end{enumerate}
$(\Omega,\sal,\nu)$ is called an internal, finitely-additive measure space.

\label{definternal}
\end{definition}

The continuity condition of Definition \ref{measure} c) is automatically satisfied  if the superstructure embedding is sufficiently saturated. 

\begin{prop} Assume that the superstructure embedding is $\kappa$-saturated. Let $\mathcal{S}_0 \subset \mathcal{S}$ be a shrinking set with empty intersection and suppose the cardinality of  $\mathcal{S}_0$ is at most $\kappa$. This is for example satisfied if $\mathcal{S}_0$ is countable. Then there is a set $A \in \mathcal{S}_0$ with $A = \varnothing$. 
\label{empty}
\end{prop}
{\em Proof.} This follows from the saturation property of the superstructure embedding $^*$. Assume that $\mathcal{S}_0$ is shrinking and all sets $A \in \mathcal{S}_0$ are non-empty. Then the intersection over all such $A$ is also non-empty, a contradiction. 
$\hfill \square$\\

This gives $\sigma$-additivity: if a countable union satisfies $\cup_{n\in\nn} A_n \in \sal$ then Proposition \ref{empty} implies that the union is finite, so that $\sigma$-additivity follows from finite additivity. But countable unions of internal sets usually fail to be internal. On the other hand, additivity holds for $^*$-finite unions, i.e. for a sequence of disjoint subsets $A_n \in \sal$ and any $N \in \nsn$, the hyperfinite union satisfies $\bigcup_{n=1}^N A_n \in \sal$ and $\nu(\bigcup_{n=1}^N A_n)= \sum_{n=1}^N \nu(A_n)$. \\

Obviously, a standard measure space $(X,\ral,\mu)$ as described in Section \ref{measure} extends to an internal, finitely additive measure space $({^*X},{^*\ral},{^*\mu})$. Another example are {\em hyperfinite} internal measure spaces which turn out to be particularly useful.

\begin{prop} Let $\Omega$ be a hyperfinite set and $\sal$ an internal and covering subring of $\mathcal{P}(\Omega)$. Then $\sal$ is the algebra  ${^*\mathcal{P}}(\Omega)$ of internal subsets of $\Omega$. Any internal function $\nu: \Omega \rightarrow {^*K}$ defines a measure on the singleton sets and can be uniquely extended to $\sal$. This defines an internal, hyperfinite measure space $(\Omega,\sal,\nu)$.
\label{hyperfinite}
\end{prop}
{\em Proof.} The corresponding statement  is obvious for any finite set $\Omega$. Then the assertion follows from the Transfer Principle.$\hfill\square$\\

\noindent {\em Examples: } Let $\Omega={^*\zp}/(p^N)$ for some prime $p \neq 2$ with infinite $N \in \nsn$
and $K=\qp$. Then any internal sequence $(\nu_i)_{0\leq i < p^N}$ with values in ${^*\qp}$ defines an internal finitely additive measure $\nu$ on $\sal={^*\mathcal{P}}(\Omega)$.

\begin{enumerate}
\item If $\nu$ is translation-invariant, then $\nu$ must be constant and if $\nu$ is also normalized, we obtain the {\em internal Haar measure} defined by $\nu_i = p^{-N}$. We have $\nu(a+p^n \Omega)=p^{-n}$ for all $n\leq N$, but since $\nu$ is not finitely bounded, it does not induce a standard measure. 
\item Set $\nu_i=1$ if $i=0,1,\dots, \frac{p^N-1}{2}$ and $\nu_i=-1$ if $i=\frac{p^N+1}{2},\dots,p^N-1$. Let $n\in \nn$ be standard. Then $\nu(a+p^n \Omega)=1$ for $a=0,1,\dots,\frac{p^n-1}{2}$ and $\nu(a+p^n \Omega)=-1$ for $a=\frac{p^n+1}{2}, \dots, p^n -1$. 
\item Set $\nu_i=(-1)^i$ for $i=0,1,\dots,p^N-1$. Then $\nu(a+p^n \Omega)=(-1)^a$  for $a=0,1,\dots,p^n-1$.
\end{enumerate} 

In the following, we reserve the letters $\ral$ and $\al$ for rings or algebras and write $\sal$ for internal rings or algebras, for example $\sal={^*\ral}$.\\

As in section \ref{measures} above, there is a Norm function $N_{\nu} : \Omega \rightarrow \nsr_{\geq 0}$. $N_{\nu}(x)$ is defined as the infimum of all $\|A\|_{\nu}$ with $x \in A \in \Omega$. We conclude by transfer from the standard case that $\|A\|_{\nu}$ can be recovered as the supremum of the set of all $N_{\nu}(x)$ with $x \in A$.\\

The measure space $\sal$ can be extended so that it includes all $\nu$-approximable sets. This extension is particularly convenient in the nonstandard setting.

\begin{definition} Let $(\Omega,\sal,\nu)$ be an internal measure space and 
$A \subset \Omega$ a (possibly external) set. $A$ is called a Loeb null set if $N_{\nu}(x) \approx 0$ for all $x \in A$. Moreover, $A$  is called Loeb-measurable if $B \in \sal$ exists such that $A \Delta B$ is a Loeb null set. The ring of Loeb-measurable sets is denoted by $\sal_L$. 
\end{definition}

Next, we derive a standard measure with values in $K$ from an internal measure space. We need to assume that $\nu$ has values in $\ns({^* K})$. If $K$ is locally compact then this is equivalent to $\nu$ having values in $\fin({^*K})$ (see Proposition \ref{pnsfin}). If $\mu$ is globally bounded then $\mu(A)=\st_K(\nu(A))$ defines a measure on $\sal$ with values in $K$ which can be extended to $\sal_L$. Furthermore, $\|A\|_{\nu}$ is a finite hyperreal number.
\\

The following Proposition shows that the Loeb construction gives a measure space in the standard sense. 

\begin{prop}  
Let $(\Omega,\sal,\nu)$ be an internal measure space and assume that $\nu$ has values in $\ns({^* K})$ and that there exists some $C \in \rr$ such that $|\mu(A)| \leq C$ for all $A \in \sal$. Then the so-called Loeb measure $\nu_L : \sal_L \rightarrow K$ is well defined: $\nu_L(A)=\st_K(\nu(B))$ where $B \in \sal$ such that $A \Delta B$ is a Loeb null set. It has the following properties:
\begin{enumerate}
\item $\nu_L(A \cup B)=\nu_L(A)+\nu_L(B)$ for disjoint sets $A, B \in \sal_L$. (Additivity)
\item For $A \in \sal_L$, $\|A\|_{\nu_L} = \sup\ \{|\nu_L(B)|\ :\ B \subset A,\ B\in \sal_L \}  \leq C$. (Boundedness)
\item For any shrinking set $\sal_{L,0} \subset \sal_L$ (i.e., $A, B \in \sal_{L,0}  \Rightarrow A \cap B \in \sal_{L,0} $) with empty intersection (i.e., $\bigcap_{A \in \sal_{L,0} } A = \varnothing$) and cardinality at most $\kappa$ one has $\displaystyle \lim_{A \in \sal_{L,0}} \| A \|_{\nu_L} = 0$. 
 (Continuity)
\end{enumerate}
We call $(\Omega,\sal_L,\nu_L)$ the Loeb measure space associated to $(\Omega,\sal,\nu)$.
\label{loebprop}
\end{prop}

\noindent {\em Proof.} Let $B, B' \in \sal$ such that $A \Delta B$ and $A \Delta B'$ are Loeb null sets. Then $N_{\nu}(x) \approx 0$ for $x \in B \Delta B'$. The additivity of $\nu$ implies $|\nu(B)-\nu(B')| \leq |\nu(B \Delta B')| \leq \| B \Delta B' \|_{\nu} \approx 0$. Hence $\nu_L$ is well defined. 
The fact that $\sal_L$ is a covering separating ring and the additivity  follow immediately from the definition. Let $A \in \sal$ and $B \subset A$. By assumption, $\nu(B)$ is nearstandard so that $|\nu(B)| \leq C$. This implies the boundedness of $\|A\|_{\nu_L}$ for $A \in \sal_L$. 

We have to show the continuity property c). Let $\sal_{L,0} \subset \sal_L$ be a shrinking set with empty intersection and cardinality at most $\kappa$. Let $\epsilon > 0$ be standard and  define $\Omega_{\epsilon} \subset \Omega$ to be the internal subset of all $y \in \Omega$ with $N_{\nu}(y) > \epsilon$.
For every $A \in \sal_{L,0}$ there exists a set $B \in \sal$ with $B \cap \Omega_{\epsilon} = A \cap \Omega_{\epsilon}$, by the definition of Loeb measurable sets. Let $\sal_0$ be the set of these $B$'s. The assumption $\bigcap_{A \in \sal_{L,0}} A = \varnothing$ implies
 \[ \bigcap_{A \in \sal_{L,0}} (A \cap \Omega_{\epsilon}) = \bigcap_{B \in \sal_0} (B \cap \Omega_{\epsilon}) =  \varnothing .  \]
 $\sal_0$ has the same cardinality as $\sal_{L,0}$. Therefore Proposition \ref{empty} gives a set $B \in \sal_0$ and a corresponding set $A \in \sal_{L,0}$ such that $B \cap \Omega_{\epsilon} = A \cap \Omega_{\epsilon} = \varnothing$. This implies $\| A \|_{\nu_L} \leq \epsilon$ which proves the continuity of the Loeb measure.$\hfill\square$

\begin{lemma} Let $(\Omega,\sal,\nu)$ be an internal measure space and $(\Omega,\sal_L,\nu_L)$ a Loeb measure space under the assumptions of the above Proposition \ref{loebprop}. 
Let $N_{\nu_L}(x)= \inf_{y \in A \in {\sal_L}} \|A\|_{\nu_L}$ be the weight function associated to to $\nu_L$. Then $\|A\|_{\nu_L} = \st_{\rr} (\| A \|_{\nu})$ for $A \in \sal$ and $ N_{\nu_L}(x)= \st_{\rr} (N_{\nu}(x))$ for all $x \in \Omega$.
\label{stloeb}
\end{lemma}
{\em Proof.} Let $B_1,B_2 \in \sal_L$. If $B_1 \Delta B_2$ is a Loeb null set, then the additivity of $\nu_L$ implies $|\nu_L(B_1)-\nu_L(B_2)| \leq \max \{ |\nu_L(B_1 \setminus B_2)| , 
|\nu_L(B_2 \setminus B_1)| \} \leq \| B_1 \Delta B_2 \|_{\nu_L} = 0$. From this and the definition of $\| \ \|$ we conclude that $\|A\|_{\nu_L} = \st_{\rr} (\| A \|_{\nu})$ for $A \in \sal$. Furthermore, it implies $N_{\nu_L}(x) \leq \st (N_{\nu}(x))$. Suppose that this is a strict inequality for $x \in \Omega$. Then there exists $A \in \sal_L$ and $B \in \sal$ such that $A \Delta B$ is a Loeb null set and  $x \in A \setminus B$. But this implies $N_{\nu}(x) \approx 0$, a contradiction.$\hfill\square$

\subsection{The standard-part map} 
For any topological space $X$ with extension ${^*X}$ there exists the standard-part map $\st_X : \ns({^*X}) \rightarrow X$ (see section \ref{nonst}) which gives a transition between nonstandard and standard spaces. We will consider the standard-part map for the real numbers, for complete non-Archimedean valued fields and for measurable spaces where  
the measurable sets form the base of a zero-dimensional Hausdorff topology.\\

First, we give some compatibility properties.

\begin{lemma} Let $(X,\ral,\mu)$ be a measure space and $({^*X},{^*\ral},{^*\mu})$ the corresponding internal measure space.
\begin{enumerate}
\item Let $y \in \ns({^* K})$. Then $|\st_K(y)| = \st_{\rr}(|y|)$.

\item For $A \in \ral$, one has $\| A\|_{\mu} = \| {^*A}\|_{^*\mu}$.

\item $N_{^*\mu} = {^*(N_{\mu})}$ and $N_{^*\mu} (y)= \inf_{y \in A \in {^* \ral}} \|A\|_{^*\mu}$.
\item Let $y \in \ns({^* X})$. Then $N_{\mu}(\st_X(y)) \geq \st_{\rr}(N_{^*\mu}(y))$.
\end{enumerate}
\label{compat}
\end{lemma}
{\em Proof. } a) Since $y \in \ns({^* K}) \subset \fin({^* K})$, $|y|$ is a finite hyperreal number and  $\st_{\rr}(|y|)$ is well defined. 
The absolute value $|\ |: K \rightarrow \rr$ is continuous which implies that the absolute values of the approximate elements $y$ and $\st_K(y)$ are also approximate (see \cite{lr} 21.11). This gives the equality.\\
b) The assertion follows from the definition of $\| \ \|$ and  a general fact on nonstandard extensions of sets which are defined by formulas (\cite{lr} 7.5): $$^*\{ |\mu(B)|: B \subset A,\ B\in \ral \}=\{ |{^*\mu}(B)|: B \subset {^*A},\ B \in {^*\ral} \}$$
c) It follows from the definition that $N_{^*\mu}$ is the extension of $N_{\mu}$.\\
d) Let $x=\st_X(y) \in A$. If $x\in A$ with $A \in \ral$, then $y \in {^* A}$ by definition of the standard-part map. Since $N_{\mu}(x)$ is the infimum of all such $\|A\|_{\mu}$, and by part b) $\|A\|_{\mu} = \|{^*A}\|_{^*\mu}$, we conclude $N_{\mu}(x) \geq \st_{\rr}(N_{^* \mu}(y))$.$\hfill\square$\\

\noindent {\em Remark:} Equality can not be deduced in d) because $y$ may be contained in additional internal measurable subsets which are not of type ${^*A}$ for $A \in \ral$. Another argument uses the fact that $N_{\mu}$ is only upper-semicontinuous (see \cite{rooij} 7.6). The elements $y$ and $x=\st(y)$ are approximate and since $N_{\mu}$ is upper-semicontinuous at $x$, only the inequality $N_{\mu}(x) \geq N_{^*\mu}(y)$ holds. \\

The following Lemma shows that $\st_X$ is defined almost everywhere.

\begin{lemma}
Let $(X,\ral,\mu)$ be a measure space, $({^*X},{^*\ral},{^*\mu})$ the corresponding internal measure space and $A \in \ral$ a measurable set. 
Then ${^* A} \setminus \ns({^*X})$ is a Loeb null set, i.e.\ $N_{^*\mu}(y) \approx 0$ for $y \in {^* A} \setminus \ns({^*X})$.
\label{nullset}
\end{lemma}
{\em Proof.} We know from Proposition \ref{compact} that $X_{\epsilon} \cap A$ is compact for standard $\epsilon > 0$. It follows that ${^* X_{\epsilon}} \cap {^* A} \subset \st^{-1}(X_{\epsilon} \cap A)$ (see Proposition \ref{topol}c) and the latter is a subset of $\ns({^*X})$. 
Using \cite{lr} 7.5  one obtains ${^* X_{\epsilon}}= {^*(X_{\epsilon})} = ({^*X})_{\epsilon}$ which is the set of all $y \in {^*X}$ with $N_{^*\mu}(y) \geq \epsilon$.
Hence $({^* X})_{\epsilon} \cap {^* A} \subset \ns({^*X})$ for all $\epsilon >0$ and thus $N_{^*\mu}(y) \approx 0$ for $y \in {^* A} \setminus \ns({^*X})$.
$\hfill\square$\\

\noindent We now give a non-Archimedean analogue of a similar Theorem on Radon spaces (see \cite{anderson82} 3.3).

\begin{theorem} Let $(X,\al,\mu)$ be a measure space, $\al$ an algebra, $(X,\almu,\mu)$ the extended standard measure space and $({^*X},{^*\al}, {^*\mu})$ the internal measure space. Assume that ${^*\mu}$ has values in $\ns({^*K})$ and let $({^*X}, {^*\al}_L, {^*\mu}_L)$ be the corresponding Loeb measure space. Then $\st_X : ({^*X},{^*\al}_L, {^*\mu}_L)  \rightarrow (X,\almu,\mu)$ is defined outside a Loeb null set. Furthermore, $\st_X$ is a measurable and measure-preserving map.
\label{stpart}
\end{theorem}
{\em Proof.} Since we assumed that $X$ is measurable, ${^*X} \setminus \ns({^*X})$ is a Loeb null set by Lemma \ref{nullset}, and $\st =\st_X$ is defined ${^*\mu}_L$-almost everywhere. ${^*\mu}$ is globally bounded by $\| X \|_{\mu}=\|{^*X}\|_{^*\mu}$.

Now let $A \in \al$, i.e.\ $A$ is clopen in the $\al$-topology. It follows from Corollary \ref{clopen} that $\st^{-1}(A) = {^*A} \cap \ns({^*X})$ so that $\st^{-1}(A)$ is Loeb measurable and 
$$\mu(A) = {^* \mu}({^* A}) =  {^* \mu}_L({^* A} \cap \ns({^*X})) = {^* \mu}_L(\st^{-1}(A))$$
A general set $A \in \almu$ can be approximated by $B_{\epsilon} \in \al$, so that for any standard $\epsilon>0$, $( A \Delta B_{\epsilon}) \subset (X \setminus X_{\epsilon})$. This yields 
$$\st^{-1}(A \Delta B_{\epsilon}) \subset \st^{-1} (X \setminus X_{\epsilon}) \subset {^* X} \setminus {^* X}_{\epsilon}$$
The latter subset relation is true since $X_{\epsilon}$ is compact (see Proposition \ref{compact}) and therefore ${^* X}_{\epsilon} \subset \st^{-1}(X_{\epsilon})$, by Proposition \ref{topol}. This implies 
 $$N_{^*\mu}(y) \leq \epsilon$$
for $ y \in \st^{-1}(A \Delta B_{\epsilon})=\st^{-1}(A) \Delta st^{-1} (B_{\epsilon}) = \st^{-1}(A) \Delta (^*B_{\epsilon} \cap \ns({^*X}))$.
We have $st^{-1} (B_{\epsilon}) ={^*B_{\epsilon}} \cap \ns({^*X})$ and $N_{^* \mu}(y) \approx 0$ for $y \in {^*X} \setminus \ns({^*X})$ which is a Loeb null set. Hence $N_{^*\mu}(y) \leq \epsilon$ for $y \in \st^{-1}(A) \Delta {^*B_{\epsilon}}$. 
 
Applying the countable comprehension principle to the internal sequence $({^*B_{\frac{1}{n}}})_{n \geq 1}$ gives a $B \in {^*\al}$ such that $N_{^*\mu}(y) \approx 0$ for $y \in (\st^{-1}(A) \Delta (B \cap \ns({^*X})))$ and hence also for 
$y \in (\st^{-1}(A) \Delta B)$.
This implies that  $\st^{-1}(A)$ is Loeb measurable  and ${^* \mu_L}(B) = {^* \mu_L}(st^{-1}(A))$. By definition of the Loeb measure, we have   ${^* \mu_L}({^* A})= {^* \mu_L}(B)$. Then the assertion follows from  
$\mu(A) = {^* \mu}({^* A})$.$\hfill\square$\\

\subsection{Liftings of Measurable Functions}

\begin{figure*}[h]
\begin{center}
\begin{tikzpicture}[ node distance=3cm, auto]

\matrix (m) [matrix of math nodes, row sep=3em, column sep=3em] 
{ \Omega & &  {^*K} \\
& & K \\ }; 
\path[->] (m-1-1) edge node[auto] {$F$ } (m-1-3);
\path[->] (m-1-1) edge node[auto] {$f$ } (m-2-3);
\path[->,dashed] (m-1-3) edge node[auto] {$\st_K$ } (m-2-3); 
\end{tikzpicture}

\caption {One-legged lifting of $f$. The standard-part map $\st_K$ is defined on the subset $\ns({^* K})$.}
\label{onelegged}  
\end{center}
\end{figure*}
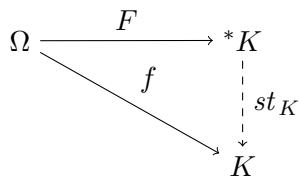

We want to show that a Loeb-measurable function $f: {\Omega} \rightarrow K$ can be lifted to an internal measurable function $F: {\Omega} \rightarrow {^*K}$ with hyperfinite image 
(see Figure \ref{onelegged}). \\
Such a function is $^*$-simple, i.e. there exist $N \in \nsn$, measurable sets $A_1,A_2,\dots,A_N$  and $y_1,y_2,\dots,y_N \in {^*K}$ with $F(x)=\sum_{i=1}^N y_i \chi_{A_i}$, where 
$\chi_{A_i}$ denotes the characteristic function of $A_i$. \\

The following Theorem shows that all Loeb-measurable functions can be lifted and that the standard part of internal measurable functions is Loeb-measurable.

\begin{theorem} Let $(\Omega,\sal,\nu)$ be an internal measure space. Assume that ${\nu}$ has values in $\ns({^*K})$ and is globally bounded. Let $({\Omega}, {\sal}_L, {\nu}_L)$ be the associated Loeb measure space.
\begin{enumerate}
\item Let $f: {\Omega} \rightarrow K$ be ${\sal}_L$-measurable and suppose that $K$ is separable.
   Then there is an internal ${\sal}$-measurable and $^*$-simple function $F:{\Omega} \rightarrow {^* K}$  such that $f(x)=\st_K(F(x))$ holds for ${\nu}_L$-almost every $x\in {\Omega}$.
\item Conversely, if $F: \Omega \rightarrow {^* K}$ is an internal $\sal$-measurable function such that $F(x) \in \ns({^* K})$ for ${\nu}_L$-almost every $x\in {^*X}$, then $f:=\st_K \circ F$ is defined $\nu_L$-almost everywhere and ${\sal}_L$-measurable.
\end{enumerate}
\label{firstlifting}
\end{theorem}

{\em Proof.} 
(a) We follow Anderson's proof for real-valued measure spaces (see \cite{anderson82} 5.3). Choose a  countable base $U_1=K, U_2, U_3 \dots$ of clopen sets in $K$, take their preimages $f^{-1}(U_n) \in \sal_L$ and replace them by $A_n \in \sal$ such that $\nu_L(f^{-1}(U_n)\Delta A_n) = 0$. Then define an approximating sequence of $\sal$-measurable simple functions $f_n : \Omega \rightarrow {^*K}$ such that $f_n(A_k) \subset {^*U_k}$ for all $k \leq n$. 
By countable comprehension, $(f_n)_{n \in \nn}$ extends to an internal sequence and $F$ can be defined as $f_{N}$ for some infinite $N\in\nsn$. $F$ is $\sal$-measurable with hyperfinite image in ${^* K}$ and one shows that $\st (F(x))=f(x)$ for $x \in \Omega \setminus \bigcup_{n=1}^{\infty}\ (f^{-1}(U_n)\Delta A_n)$. Although $\sal_L$ is not a $\sigma$-algebra, the countable union of Loeb null sets is again a null set. \\
(b) Let $A \subset K$ be a clopen set. Then 
$$f^{-1}(A)=F^{-1}(\st^{-1}(A)) = F^{-1}({^*A} \cap \ns({^*K})) = F^{-1}({^* A}) \cap F^{-1}(\ns({^*K})) . $$
Since $F$ is $\sal$-measurable, $F^{-1}({^* A})\in \sal$. The assumption on $F$ ensures that the intersection with $F^{-1}(\ns({^*K}))$ reduces this set only by a $\nu_L$-null set and hence $f^{-1}(A) \in \sal_L$.  \hspace*{1cm} $\hfill\square$\\

Besides liftings of Loeb-measurable functions on an internal measure space $\Omega$, there are also liftings of {\em standard} measurable functions $f:X \rightarrow K$. Natural candidates are the nonstandard extension ${^* f}$ and a lifting $F$ of the composition $f \circ \st$ as constructed above (see  
Figure \ref{two-legged}).

\begin{figure*}[h]
\begin{center}
\begin{tikzpicture}[ node distance=3cm, auto]

\matrix (m) [matrix of math nodes, row sep=3em, column sep=3em] 
{ {^*X} & & {^*K} \\
X & & K \\ }; 
\path[->] (m-1-1) edge node[auto] {${^*f}, F$ } (m-1-3);
\path[->,dashed] (m-1-1) edge node[auto] {$\st_X$} (m-2-1);
\path[->]  (m-2-1) edge node[auto] {$f$ } (m-2-3);
\path[->,dashed] (m-1-3) edge node[auto] {$\st_K$ } (m-2-3); 
\end{tikzpicture}

\caption {Two-legged lifting of $f$. The standard-part maps $\st$ are defined on the subsets $\ns({^* X})$ resp.\ $\ns({^* K})$.}
\label{two-legged}  
\end{center}
\end{figure*}

\begin{theorem} Let $(X,\al,\mu)$ be a measure space, $\al$ an algebra, $(X,\almu,\mu)$ the extended measure space, $({^*X},{^*\al},{^*\mu})$ the internal measure space. Assume that  ${^*\mu}$ has values in $\ns({^*K})$. Let $({^*X}, {^*\al}_L, {^*\mu}_L)$ be the corresponding Loeb measure space. Let $f: X \rightarrow K$ be a $\almu$-measurable function. Then ${^* f}: {^*X} \rightarrow {^*K}$ is ${^*(\almu)}$- and ${^*\al}_L$- measurable and 
$\st_K({^* f}(y))=f(\st_X(y))$ for ${^* \mu}_L$-almost all $y \in {^*X}$, i.e.\ the diagram in Figure \ref{two-legged} commutes almost everywhere.
\label{commute}

\end{theorem}
{\em Proof. } By transfer,  ${^* f}$ is $({^*\almu}, {^*B(K)})$-measurable. Since ${^*(\almu)} \subset  {^*\al}_L$, ${^*f}$ is also measurable with respect to ${^*\al}_L$.

Since $f$ is $\almu$-continuous (see Corallary \ref{cont}), the proof of the commutative diagram is easier than in the real case (see \cite{anderson82} 3.7). Let $y \in \ns({^* X})$ and set $x=\st_X(y)$. Then $y \approx x$ and the continuity implies ${^*f}(y) \approx f(x)$ so that $\st_K({^*f}(y))=f(x)=f(\st_X(y))$. ${^*X} \setminus \ns({^* X})$ is a null set (cf. Proposition \ref{nullset}) and this gives the assertion. $\hfill\square$\\

{\em Remark.} The composition map $f \circ \st_X : {^* X} \rightarrow K$ is ${^* \al}_L$-measurable and has a ${^* \al}$-measurable lifting $F:{^* X} \rightarrow {^* K}$ (see Theorem \ref{firstlifting}). What are the main differences between ${^* f}$ and $F$ ? Of course, ${^* f}$ is the natural nonstandard extension of $f$ and 
${^* \al}_L$-measurable. In contrast, $F$ is obtained by a non-canonical construction from $f$, but it has additional favourable properties: $F$ is ${^* \al}$-measurable and even ${^*}$-simple.

\subsection{Hyperfinite Spaces}

A key technique in nonstandard analysis is {\em hyperfinite approximation}. A well-known example from nonstandard real analysis is the representation of the compact interval $X=[0,1]$ by the hyperfinite set $Y=\{0,\frac{1}{N},\frac{2}{N}, \dots, \frac{N-1}{N},1\}$,  where $N \in \nsn$ is  some infinite natural number. There is a  hyperfinite measure space defined on $Y$ and the standard-part map is measurable and measure-preserving. Standard Lebesgue integrals on $X$ can be obtained by a hyperfinite summation over $Y$.  This was generalized to Radon probability spaces (see \cite{anderson82}, \cite{albeverio1986}) and it can be shown that Radon measures are push-downs of hyperfinite measures spaces. We show a corresponding statement for measures with values in a complete non-Archimedean field $K$.\\

For a measure space $(Y,\mathcal{S},\tau)$ and a map $p: Y \rightarrow X$, the {\em push down} of $Y$ to $X$ via $p$ induces a measure space on $X$. A set $A \subset X$ is measurable if $p^{-1}(A) \in \mathcal{S}$ and the measure $p(\tau)$ is defined by $p(\tau)(A)=\tau(p^{-1}(A))$. 

\begin{theorem} Let $\al$ be an algebra, $(X,\al,\mu)$ a measure space,  $(X,\almu,\mu)$ the extended measure space and 
$({^*X},{^*\al},{^*\mu})$ the corresponding internal measure space. Assume that ${^*\mu}$ has values in $\ns({^*K})$. Then there is a hyperfinite partition of ${^*X}$ such that all partition classes are contained in ${^*\al}$, an equivalence relation $\thicksim$ on ${^* X}$ defined by the partition and a hyperfinite set $Y={^* X}/\thicksim$. It induces a measure $\nu={^*\mu}/\thicksim$ and an internal hyperfinite measure space $(Y,\sal,\nu)$ such that all internal subsets of $Y$ are measurable. Let $(Y,\sal_L,\nu_L)$ be the corresponding Loeb measure space. Then the standard-part map $\st_Y: Y  \rightarrow X$ is well-defined $\sal_L$-almost everywhere and
$(X,\almu,\mu)$ is  the push-down of $Y$ to $X$ via $\st_Y$. In particular, $\st_Y$ is measurable and measure-preserving.
\label{hyperfinite-rep}
\end{theorem}

{\em Proof.} Let $A_1, A_2, \dots, A_n$ be clopen sets in $X$. As in \cite{albeverio1986} 3.4.10, let $P_{A_1,A_2,\dots A_n}$ be the set of hyperfinite partitions of ${^* X}$ into ${^* \mu}$-measurable sets such that each ${^* A_i}$ ($i=1,\dots,n$) is a disjoint union of partition classes. If the nonstandard extension is sufficiently saturated (the usual countable saturation suffices for separable spaces $X$), we obtain a hyperfinite partition $P=\{R_1, R_2, \dots, R_N \}$ of  ${^* X}$ into ${^* \mu}$-measurable sets such that each set ${^* A}$ (where $A$ is clopen in $X$) is a disjoint union of $R_i$'s. $P$ defines an equivalence relation $\thicksim$ on ${^* X}$ and $Y:={^* X}/\thicksim$ is a hyperfinite set with $N$ elements. Let $\sal \subset {^*\al}$ be the internal set of all hyperfinite unions of partition classes $R_i$ and for $B \in \sal$ define $\nu(B)={^*\mu}(B)$. Then $(Y,\sal,\nu)$ is an internal hyperfinite measure space where all internal subsets of $Y$ (resp.\ their corresponding union of partition classes in ${^*X}$) are measurable. 

We show that the standard-part map $\st=\st_X: \ns({^* X}) \rightarrow X$ factorizes via the equivalence relation $\thicksim$ which defines $Y$. Anderson \cite{anderson82} calls this property {\em S-separating} and it means that $\st(y)$ does not depend on the choice of $y \in R_i$. To this end, suppose that $y,y' \in R_i$. For all clopen sets $A \subset X$, the elements $y$ and $y'$ belong to the same sets ${^* A}$, depending on whether $R_i \subset  {^* A}$ holds or not. This implies that they are in the same monad and $\st(y)=\st(y')$, if $y \in \ns({^* X})$, or equivalently, $y' \in \ns({^* X})$.

Let $(Y,\sal_L,\nu_L)$ be the Loeb measure space which corresponds to $(Y,\sal,\nu)$. Recall that $A \in \sal_L$ if $A$ is a (possible external) union of $R_i$'s and $B \in \sal$ exists with $N_{\nu}(y) \approx 0$ for 
$y \in  A \Delta B$. In that case, $\nu_L(A)=\st_K(\nu(B)) \in K$. Subsequently, we prove that $\st_Y$ is $\sal_L$-measurable and measure-preserving. First, let $A \in \al$ be a clopen set. Then $\st_Y^{-1}(A)= {^* A} \cap \ns({^* X})$ which is a (possibly external) union of partition classes. Since ${^*X} \setminus \ns({^* X})$ is a ${^* \mu}$- and a $\nu_L$- null set, we obtain  $\st_Y^{-1}(A) \in \sal_L$ and also $\nu_L(\st_Y^{-1}(A)) = {^*\mu}_L({^* A} \cap \ns({^* X})) = {^* \mu}({^* A})=\mu(A)$. For general $A \in \al_{\mu}$, one proceeds as in the proof of Theorem \ref{stpart}.

Finally, we prove that $A \subset X$ is $\almu$-measurable if $\st^{-1}(A) \in \sal_L$. Let $\epsilon > 0$ be any standard real number and $X_{\epsilon}$ the set of all $x \in X$ with $N_{\mu}(x) > \epsilon$. It suffices to find a set $C \in \almu$ such that $A \cap X_{\epsilon} = C \cap X_{\epsilon} $ (see \cite{rooij} 7.3 and 7.8). 
The assumption yields a $B \in \sal$ such that $ B \Delta \st^{-1}(A) $ is a Loeb null set. 
Hence $N_{^*\mu}(x) \leq \epsilon$ for $ x \in B \Delta \st^{-1}(A)$ so that the intersection of this set with ${^* X}_{\epsilon}$ is empty and $\st^{-1}(A) \cap  {^* X}_{\epsilon} = B \cap {^* X}_{\epsilon}$.  We apply the $\st$-map and get $A \cap X_{\epsilon} = \st(B) \cap X_{\epsilon}$ since $\st({^*X_{\epsilon}})= X_{\epsilon}$ by Propositions \ref{compact} and \ref{topol}. 
Since $B$ is internal, $\st(B)$ is closed (see \cite{lr} 28.7) and as $X_{\epsilon}$ is compact, we obtain that $A \cap X_{\epsilon}$ is closed. The same argument applies to $X \setminus A$ since ${^* X} \setminus \st^{-1}(A) \in \sal_L$. This implies that $(X \setminus A) \cap X_{\epsilon}$ is closed. 
By the above, $X \setminus (A \cap X_{\epsilon})$ is $\almu$-open and hence a union of $\almu$-clopen sets $A_i$. 
The intersection  $(X \setminus (A \cap X_{\epsilon})) \cap X_{\epsilon} =  (X \setminus A) \cap X_{\epsilon}$ is not only closed, but also compact (since $X_{\epsilon}$ is compact). 
It can hence be covered by a finite union of $A_i$'s. The complement of this union in $X$, which we call $C$, is a clopen set and has the desired property $A \cap X_{\epsilon} = C \cap X_{\epsilon} $. $\hfill\square$ \\

{\em Remark:} For an ultrametric space $X$  with nonstandard extension ${^*X}$, the hyperfinite space $Y$ can be chosen as the set of all balls of radius $\epsilon$ in ${^*X}$, where $\epsilon>0$ is a  fixed infinitesimal number. The balls form a partition of ${^*X}$ which is $S$-separating since each ball $B_{\epsilon}(y)$ is a subset of the monad of $\st_X(y)$ for $y \in \ns({^*X})$. The above Theorem says that a measure on an ultrametric space $X$ is determined by a hyperdiscrete measure on the hyperfinite space $Y$ of infinitesimal balls. 

\section{Integration}
\label{sec:int}

With the preceding results, it is not surprising that nonstandard extensions can also be used for the integration of functions with values in non-Archimedean fields. 
For the convenience of the reader we recall some facts from the standard theory \cite{rooij}. \\

For a measure space $(X,\ral,\mu)$ and a field $K$ as above, the seminorm of a function $f: X \rightarrow K$  is defined as $\| f \|_{\mu} = \sup_{x\in X} |f(x)| \cdot N_{\mu}(x) \in [0,\infty[ \cup \{\infty\}$. Let $\chi_A : X \rightarrow K$ denote the characteristic function of $A \subset X$. For a measurable set $A \in \ral$, we have $\|\chi_A\|_{\mu} = \| A \|_{\mu}$. $f$ is a simple function (step function), if $A_1,A_2,\dots,A_n \in \ral$ and $x_1,x_2,\dots,x_n \in K$ exist with $f(x)=\sum_{i=1}^n x_i \chi_{A_i}$. The integral is a functional on the space of simple functions defined by $\int_X f(x)\ d\mu = \sum_{i=1}^n x_i \mu(A_i)$ and it satisfies the inequality $| \int_X f(x)\ d\mu | \leq \| f \|_{\mu}$. The space of $\ral$-simple functions can be completed w.r.t.\ $\|\ \|_{\mu}$ and a function $f: X\rightarrow K$ 
is called $\mu$-integrable, if a sequence $(f_n)_{n\in \nn}$ of $\ral$-simple functions exists such that $\lim_{n \rightarrow \infty} \| f - f_n \|_{\mu} = 0$. The integral $\int_X f(x)\ d\mu$ is defined as $\lim_{n\rightarrow \infty} \int_X f_n(x)\ d\mu$. One can show that $\chi_A$ is $\mu$-integrable if and only if $A \in \ralmu$ and this also provides an alternative way to define the extended ring $\ralmu$. For simple functions the notions of $\mu$-{\em integrability} and $\ralmu$-measurability coincide, but for general functions integrability requires an additional boundedness condition as the following Theorem (\cite{rooij} 7.12 and Corollary \ref{cont}) shows:

\begin{theorem} Let ($X,\ral,\mu)$ be a measure space and $K$ a field as above. A function $f: X\rightarrow K$ is $\mu$-integrable if and only if $f$ is locally $\ralmu$-measurable and for every $\epsilon>0$, the set $\{ x\in X : |f(x)|N_{\mu}(x) \geq \epsilon \}$ is $\ralmu$-compact, hence contained in some set $X_{\delta}$ for $\delta >0$.
\label{integral}
\end{theorem}

Now we study the corresponding nonstandard representation of integrals and show that integrals of $^*$-simple are sufficient.

Let $(\Omega,\sal,\nu)$ be an internal measure space, $N \in \nsn$ and $A_1,\dots,A_N \in \sal$.
Then the $^*$-simple function $F=\sum_{i=1}^N y_i \chi_{A_i}$ is $^*$-integrable and $\int_{\Omega} F(x)\ d\nu = \sum_{i=1}^N y_i \mu(A_i) \in {^*K}$. The seminorm on standard functions extends to internal functions $F: \Omega \rightarrow {^* K}$. For $^*$-simple functions $F$ we have $\| F\|_{\nu} \in \nsr_{\geq 0}$ (not $\infty$, but not necessarily finite). It is useful to put the following restriction on ${^*}$-integrability in order to relate it to standard integrals. 

\begin{definition} Let $(\Omega,\sal,\nu)$ be an internal measure space and $F: \Omega \rightarrow {^*K}$ a $^*$-simple function. $F$ is called $S$-integrable if the following conditions are satisfied:
\begin{enumerate}
\item $\|F\|_{\nu}$ is a finite hyperreal number, and 
\item If $A \in \sal$ with $\|A\|_{\nu} \approx 0$ then $\| F \cdot \chi_{A} \|_{\nu} \approx 0$, and
\item If $A \in \sal$ with $F \cdot \chi_A \approx 0$ then $\| F \cdot \chi_{A} \|_{\nu} \approx 0$.
\end{enumerate}
\end{definition}

We have adopted the definition of $S$-integrability \cite{anderson76} from real nonstandard measure theory. If $\sal$ is an algebra such that $\|\Omega\|_{\nu}$ is finite, then c) is automatically satisfied. We remark that a $^*$-simple function $F=\sum_{i=1}^N y_i \chi_{A_i}$ is $S$-integrable if $|y_i|$ and $\|A_i\|_{\nu}$ are finite for all $i=1,\dots,n$.
But for example $\delta$-functions (which exist as nonstandard {\em functions}) are not $S$-integrable. \\
We already know from Theorem \ref{firstlifting} that a Loeb measurable function can be lifted to a simple function, and conversely, the standard part of a simple function is Loeb measurable. The following Theorem relates $S$-integrability to Loeb integrability. Since Loeb measurable functions have values in a standard non-Archimedean field, the conventional  definition of integrals can be used here.

\begin{theorem} Let $(\Omega,\sal,\nu)$ be an internal measure space, assume that $\nu$ has values in $\ns({^*K})$ and let $(\Omega,\sal_L,\nu_L)$ be the corresponding Loeb space. 

\begin{enumerate}
\item Let $f: \Omega \rightarrow K$ be integrable (in the conventional sense w.r.t.\ the Loeb measure). Then $f$ has a $S$-integrable and $^*$-simple lifting $F : \Omega \rightarrow {^*K}$.

\item Let $K$ be a locally compact and $F: \Omega \rightarrow {^* K}$ a $S$-integrable and $^*$-simple function such that $f= \st_K \circ F$ is defined outside a Loeb null set. Then $f$ is integrable w.r.t.\ the Lob measure.
\end{enumerate}
Under the hypotheses of a) or b),
\[ \int_{\Omega} f(x)\ d\nu_L = \st_K \int_{\Omega} F(x)\ d\nu \]
\label{internal-int}
\end{theorem}

{\em Proof.} a): Let
$(f_n)_{n\in\nn}$ be a sequence of $\sal_L$-simple functions which converges to $f$ w.r.t.\ $\|\ \|_{\nu_L}$. It follows from the definition of $\sal_L$ that there is a $\sal$-simple function $g_n: \Omega \rightarrow K$ for each $n$ such that $f_n = g_n$ outside a Loeb null set. This gives $\| f_n - g_n \|_{\nu_L} = 0$ and $\| f_n - g_n \|_{\nu} \approx 0$. Thus
one obtains $\|f-f_n\|_{\nu_L} \approx  \|f-f_n\|_{\nu} \approx \|f-g_n\|_{\nu}$.
 
We apply countable saturation to the internal sets of $^*$-simple functions $g$ with $\|f-g\|_{\nu}<\epsilon$. Hence there exists a  $^*$-simple function $F: \Omega \rightarrow {^*K}$ such that $\|f-F\|_{\nu} \approx 0$. This implies $|f(x)-F(x)| \approx 0$ and thus $\st_K(F(x))=f(x)$ outside a Loeb null set. 
This gives a $^*$-simple lifting of $f$. 
Since $\| F\|_{\nu} \leq \max\{ \|F-f\|_{\nu}, \|f\|_{\nu} \}$, it follows that $\| F\|_{\nu}$ is finite. If $A \in \sal$ with $\| A\|_{\nu} \approx 0$, then $\| F \chi_A\|_{\nu} \leq \max \{ \|(F-f)\chi_A\|_{\nu}, \|f \chi_A\|_{\nu} \}$, which is infinitesimal. 
If $F \cdot \chi_A \approx 0$ then $\| F \cdot \chi_A \|_{\nu} \approx 0$ since $\|A\|_{\nu}$ is finite.

This shows that $F$ is $S$-integrable.\\

b): Let $F: \Omega \rightarrow {^* K}$ be a $S$-integrable lifting of $f$ and $F=\sum_{i=1}^N y_i \chi_{A_i}$. Since $\nu$ is an internal measure with values in $\ns({^*K})$,  all $\|A_i\|_{\nu}$ are finite and there is a standard real number $M>0$ such that $N_{\nu}(x) < M$ for all $x \in \bigcup_{i=1}^N A_i$. Let $\epsilon>0$ be any standard real number. We claim that there exists a $\sal$-simple function $g :\Omega \rightarrow K$ such that $\| F - g\|_{\nu} < \epsilon$ which implies $\| f - g\|_{\nu_L} \leq \epsilon$ and proves the assertion.\\

First, assume that for all $D \in \rr$ there exists $x \in \Omega$ such that $|F(x)| > D$ and $|F(x)| N_{\nu}(x) > \epsilon$. Since $F$ is $S$-integrable, the seminorm $C=\| F\|_{\nu}$ is finite. This implies $N_{\nu}(x) < \frac{C}{D}$ which gives a set $A \in \sal$ with $x \in A$, $\| A\|_{\nu} < \frac{C}{D}$ and $\| F \chi_A\|_{\nu} > \epsilon$. Then the countable saturation principle yields the existence of  a set $A \in \sal$ with $\| A\|_{\nu} \approx 0$ and $\| F \chi_A\|_{\nu} > \epsilon$ which contradicts the second property of $S$-integrability. Hence there is a $D \in \rr$ such that for all $x \in \Omega$ one has $|F(x)| \leq D$ or $|F(x)| N_{\nu}(x) < \epsilon$. We may therefore take an approximating simple function $f_n$ which is zero on $x \in \Omega$ with $|F(x)| > D$. \\
Since $K$ is locally compact, the ball with radius $D$ and center $0$ is the finite disjoint union of balls $B_i$ of radius $r$ with $0<r<\frac{\epsilon}{M}$. We define $g$ to be constant on the $\sal$ measurable sets $F^{-1}(B_i)$ and zero elsewhere.  The value of $g$ on $F^{-1}(B_i)$ is an arbitrary element of $B_i$. $g$ is a $\sal$-simple function. Then $|F(x)-g(x)| < \frac{\epsilon}{M}$ and hence $|F(x)-g(x)| N_{\nu}(x) < \epsilon$ holds for all $x$ with $|F(x)| \leq B$. This implies $\| F - g\|_{\nu} < \epsilon$ as claimed.\\

\noindent It remains to prove the integral formula. Under the hypothesis of a) or b), we have a $^*$-simple lift and S-integrable lifting $F: \Omega \rightarrow {^*K}$ of the integrable function $f: \Omega \rightarrow K$.   Moreover, for any standard $\epsilon > 0$ we have a $\sal$-simple function $g : \Omega \rightarrow K$ with $\| F- g\|_{\nu} < \epsilon$ and $|\int_{\Omega} g\ d\nu -\int_{\Omega} f\ d\nu_L | < \epsilon$. These two inequalities imply $|\int_{\Omega} F\ d\nu -\int_{\Omega} f\ d\nu_L | < \epsilon$ which  proves the integral formula.
$\hfill \square$ \\

Finally, we consider {\em standard } measure spaces $(X,\almu,\mu)$ with an algebra $\al$ and integrable functions $f: X \rightarrow K$. By Theorem \ref{hyperfinite-rep} we know that $X$ is the push-down of a hyperfinite measure space $Y$ via the standard-part map. The following  result shows that the integral on $X$ can be represented by a hyperfinite summation on the space $Y$. 

\begin{theorem} Let $(X,\almu,\mu)$ be the push-down of a hyperfinite measure space $(Y,\sal_L,\nu_L)$ via the standard-part map $\st_Y$ as in Theorem \ref{hyperfinite-rep}. Let $f: X \rightarrow K$ be $\mu$-integrable. Then $f \circ \st_Y$ is Loeb integrable and there is an S-integrable lift $F: Y \rightarrow {^* K}$ such that $f \circ \st_Y = \st_K \circ F$ holds $\sal_L$-almost everywhere. The following integrals coincide:

\[ \int_X f(x)\ d\mu = \int_{Y} f(\st_Y(y))\ d\nu_L = \int_Y \st_K(F(y))\ d\nu_L   = \st_K \int_{Y} F(y)\ d\nu \]
\label{intformel}
\end{theorem}

{\em Proof.} Since $f$ is $\almu$-measurable by assumption and $\st_Y$  is measurable by Theorem \ref{hyperfinite-rep}, we obtain the Loeb measurability of  $f \circ \st_Y$. We have show that $f \circ \st_Y$ is integrable and that the integrals of $f \circ \st_Y$ and $f$ coincide. The remaining statements then follow from Theorems \ref{internal-int} and \ref{commute}.  \\

We have to prove a {\em change of variable} statement.
Since $f$ is $\mu$-integrable there is a sequence of $\al$-simple functions $f_n : X \rightarrow K$ such that $\lim_{n\rightarrow \infty} \| f - f_n \|_{\mu} = 0$. Since $\st_Y$ is measurable and $\mu(A) = \nu_L(\st_Y^{-1}(A))$ for $A \in \al$, it is obvious that $f_n \circ \st_Y$ is $\sal_L$-simple and $\int_X f_n(x)\ d\mu = \int_Y f_n(\st_Y(y))\ d\nu_L$. We have to show that the sequence $f_n \circ \st_Y$ converges to $f$, i.e.\  
$\lim_{n\rightarrow \infty} \| (f - f_n)\circ \st_Y \|_{\nu_L} = 0$. It is sufficient to prove that $\| g\circ \st_Y \|_{\nu_L} \leq \| g \|_{\mu}$ for any function $g: X \rightarrow K$. To this end, we let $y \in Y$ and set $x=\st_Y(y)$. By definition of the seminorm, the above inequality reduces to $N_{\nu_L}(y) \leq N_{\mu}(x)$. For any $A \in \al$ with $x\in A$, one has $y \in \st^{-1}(A) =  {^* A} \cap \ns({^*X})$ by Corollary \ref{clopen}. $X \setminus \ns({^*X})$ is a Loeb null set and Lemma \ref{compat} yields
$\| A\|_{\mu} = \|{^*A}\|_{^*\mu} \approx  \|\st_Y^{-1}(A)\|_{\nu_L}$.  The definition of $N_{\mu}(x)$ resp.\ $N_{\nu_L}(y)$ as infimum of seminorm values of measurable sets finally shows the inequality. 
Note  that $\sal_L$ usually contains more (and in particular smaller) measurable sets than those of type $\st_Y^{-1}(A)$ so that only the desired inequality can be obtained by this argument. $\hspace*{1cm} \hfill\square$

\section{Examples and Applications}
In this section, we apply our results to the compact ultrametric spaces $X=\zp$ and $X^{\times}=\zp^{\times}$. The algebras $\al=B(X)$ and $\al=B(X^{\times})$ of clopen subsets consist of finite unions of {balls} $B_{p^{-n}}(a)=a+p^n \zp$.  \\

\subsection{$q$-adic Haar Measure}

We begin with measures with values in a $q$-adic number field $K=\qqq$ for a prime number $q \neq p$. 
There is the translation invariant {\em Haar measure} on $X$ defined by $\mu(B_{p^{-n}}(a)) = p^{-n}$. Since $|p|=1$ for the $q$-adic absolute value, we have $|\mu(A)|=1$ for any non-empty $A \in \al$. Hence $N_{\mu}(x)=1$ for all $x\in X$ and therefore $\al=\al_{\mu}$.  The following internal measure space $(Y,\sal,\nu)$ is a hyperfinite representation of $(X,\al,\mu)$ (see Theorem \ref{hyperfinite-rep}): choose any infinite $N \in \nsn$ and set $Y={^*\zp}/(p^N)$. Let $\sal$ be the algebra of all internal subsets of $Y$ 
and define $\nu$ as a normalized counting measure: $\nu(\{y \}) = p^{-N}$. Then $\nu(A)=\frac{\# A}{p^N}$ where $\# A \in \nsn$ denotes the number of elements which is a well-defined by the transfer principle. The standard-part map $\st_Y : {^*\zp}/(p^N) \rightarrow \zp$ is measurable and measure-preserving. Furthermore, $A \subset \zp$ is measurable if and only if $\st_Y^{-1}(A)$ is an internal subset of the hyperfinite space $Y$.  \\

Since $X$ is a compact space without $\mu$-null sets (other than $\varnothing$), the following conditions   are equivalent for a function $f: X \rightarrow K$: a) $f$ is $\mu$-measurable, b) $f$ is $B(X)$-continuous, and c) $f$ is $\mu$-integrable. If this is satisfied, $f$ can be lifted to an internal simple $S$-integrable function $F: Y  \rightarrow {^*K}$ where we can take $F(x+p^N{^*\zp})={^* f}(x)$. This does not depend on the choice of the representative $x$, since $f$ is continuous. Theorem \ref{intformel} implies
\begin{equation} \int_X f(x)\ d\mu = \st \left( \int_Y F(y)\ d\nu \right) = \st \left( \frac{1}{p^N} \sum_{x=0}^{p^N-1} {^*f}(x) \right)
\label{q-adic}
\end{equation}

We remark that a similar result is well known for real-valued measures \cite{cutland2000loeb}. Consider for example the compact real interval $X=[0,1]$. The corresponding Lebesgue measure space can be represented by the {\em hyperfinite time line} $Y=\{0,\frac{1}{N},\frac{2}{N}, \dots , \frac{N-1}{N} \}$, the algebra of internal subsets of $Y$ and the normalized counting measure defined by $\nu(\{y\})=\frac{1}{N}$. Any Lebesgue-integrable function $f: X \rightarrow \rr$ (or $\cc$) can be lifted to an internal simple $S$-integrable function $F: Y  \rightarrow {\nsr}$ (resp.\ ${^*\cc})$. For continuous (or more general Riemann-integrable) functions $f$, one can take $F(y)={^* f}(y)$. The following integrals coincide:
\begin{equation}
\int_X f(x)\ d\lambda = \st \left( \int_Y F(y)\ d\nu \right) = \st \left( \frac{1}{N} \sum_{k=0}^{N-1} {^*f}(\frac{k}{N}) \right)
\label{real}
\end{equation}

Although \eqref{q-adic} and \eqref{real} look very similar, there are important differences: $|\frac{1}{p^N}| = 1$ for the $q$-adic norm, whereas $|\frac{1}{N}| \approx 0$ for the real norm. In the $q$-adic case, only the internal subsets of $Y$ are measurable. For real-valued measures there are additional external Loeb measurable subsets. Needless to say that the space of Lebesgue integrable functions contains discontinuous functions which makes the construction of a $^*$-simple lifting less obvious.

\subsection{Kubota-Leopoldt $p$-adic Zeta Function}
\label{klzeta}
Now let $p\neq 2$ and  $K=\qp$. 
First, it is well known that there exists no $p$-adic Haar {\em measure} and only a Haar {\em distribution} on $X=\zp$ and $X^{\times}=\zp^{\times}$.
But for any infinite $N \in \nsn$ there is an {\em internal Haar measure} on the hyperfinite set $Y=\{0,1,\dots,p^N-1\}$ given by $\nu(\{a\})=\frac{1}{p^N}$. The integral of any internal function $F:Y \rightarrow {^*K}$ is well defined. Any standard function $f: \nn \rightarrow K$ can be uniquely lifted to $Y$ and by restriction to $\nn$, a function $f: \zp \rightarrow K$ can also be lifted to an internal function $F:Y \rightarrow {^*K}$. We call this the {\em interpolation lifting} of $f$ and an explicit representation is given by the $^*$-finite Mahler polynomial

$$ F(y) = \sum_{n=0}^{p^N-1} a_n \binom{y}{n}  \text{ where } a_n = \sum_{i=0}^n (-1)^{n-i} \binom{n}{i} f(i) $$ 

 If $f$ is continuous then 
\begin{equation}
\st_{K}(F(y)) \approx f(\st_X(y))
\label{Fcompat}
\end{equation}
Observe that there is a unique interpolation lifting $F$, but there are also other liftings satisfying the relation \eqref{Fcompat}. If the integral $\int_Y F\ d\nu$  is nearstandard then the standard part gives the {\em Volkenborn integral} \cite{volkenborn1}.  It follows from Lemma \ref{convergence} that $\st_K(\int_Y F d\nu)$ does not depend on the infinite number $N$ and the hyperfinite space $Y$.

We consider the internal function $F: Y \rightarrow {^*K}$ given by $F(y)=y^n$ where $n \in \nn$. 
By transfer of a classical formula on sums of powers and Bernoulli numbers, the integral of $F$ can be easily computed:

\[ \int_Y y^n\ d\nu = \frac{1}{p^N} \sum_{m=0}^{p^N-1} m^n =   \frac{1}{n+1} \sum_{j=0}^n \binom{n+1}{j} B_j p^{N\cdot(n+1-j)} \frac{1}{p^N} \]
The last sum runs through a finite set 
and the standard-part of all summands with $j<n$ vanishes since $\st(p^N)=0$. The remaining term equals the Bernoulli number $B_n$ which is closely related to the value of the Riemann zeta function at the negative integer $1-n$ (see \cite{leopoldt}): 
\[ \st \left(\int_Y y^n\ d\nu\right) = B_n = (-1)^{n-1} n \zeta(1-n) \]

It is well known that $B_n$ and $\zeta(1-n)$ vanish for odd $n\geq 3$.\\

Now let $Y^{\times} \subset Y$ be the subset of integers prime to $p$. We restrict $\nu$ to $Y^{\times}$ and obtain again a hyperfinite measure space. The obvious decomposition
\[ \sum_{m=0}^{p^N -1} m^n = p^n
 \sum_{m=1}^{p^{N-1}-1} m^n + \sum_{\substack{m=0 \\ p\, \nmid\, m}}^{p^N-1} m^n\]
yields the integral formula
\[ (1-p^{n-1}) \st \left(\int_Y y^n\ d\nu \right) = \st \left(\int_{Y^{\times}} y^n\ d\nu\right) = (-1)^{n-1} n (1-p^{n-1}) \zeta(1-n) \]

Below, we will use similar integrals for nonstandard versions of $p$-adic zeta functions. \\

From now on we assume that $K=\cc_p$. Elements in $\zp^{\times}$  can be uniquely written as $x=\omega(x) \langle x \rangle$ where $\omega$ is the Teichmüller character, $\omega(x)$ a $(p-1)$-st root of unity and  $\langle x \rangle \in 1+p\zp$. 
It is well known that ${\langle x \rangle}^s=\exp_p( s \log_p \langle x \rangle) $ is well-defined, continuous and analytic in $s \in \cc_p$ if $|s|<p^{(p-2)/(p-1)}$:

\begin{equation} {\langle x \rangle}^s = \sum_{n=0}^{\infty} (\log_p \langle x \rangle)^n \frac{s^n}{n!} = \sum_{n=0}^{\infty} \binom{s}{n} (\langle x \rangle - 1)^n
\label{power-exp}
\end{equation}

Let $i \in \zz/(p-1)\zz$. Then $\omega^{1-i}$ is a well-defined power of the Teichmüller character. The $p$-adic zeta functions $\zeta_{p,i}(s) = L_p(s,\omega^{1-i})$ can be defined as a $p$-adic Mellin transform using the Volkenborn integral:

\begin{equation} \zeta_{p,i}(s) = \frac{(-1)^{i-1}}{s-1}  \int_{X^{\times}} \omega(x)^{1-i} \langle x \rangle^{1-s}\ dx  
\label{kubota-org}
\end{equation}

This is basically the original definition of Kubota and Leopoldt.
 $\zeta_{p,i}(s)$ is a $p$-adic meromorphic function and $\zeta_{p,1}(s)=L_p(s,1)$ has a simple pole in $s=1$. \\
 
 Now we give nonstandard formulas for  $\zeta_{p,i}(s)$. 
\begin{align} \zeta_{p,i}(s) & = \frac{(-1)^{i-1}}{s-1}  \st \left ( \int_{Y^{\times}} \omega(y)^{1-i} \langle y \rangle^{1-s} \ d\nu \right)  \label{kubota} \\
& = \frac{(-1)^{i-1}}{s-1}  \st \left(  \int_{Y^{\times}} \omega(y)^{1-i} \sum_{n=0}^{M} (\log_p \langle y \rangle)^n   \frac{(1-s)^n}{n!}  \ d\nu \right) \label{kubota-logp} \\
& = \frac{(-1)^{i-1}}{s-1}  \st \left( \int_{Y^{\times}} \omega(y)^{1-i} \sum_{n=0}^M \binom{1-s}{n} (\langle y \rangle -1)^{n} \ d\nu \right) \label{kubota-exp}  
\end{align}
 
One has $|\log_p\langle y \rangle| \leq \frac{1}{p}$ and $ |\langle y \rangle -1| \leq \frac{1}{p}$. Furthermore,  $|\frac{(1-s)^n}{n!}| \leq p^{n(p-2)/(p-1)} | \frac{1}{n!}|$ and 
$|\binom{1-s}{n}| \leq p^{n(p-2)/(p-1)} |\frac{1}{n!}|$. An easy computation shows that there exists a $M \in \nsn$ such that for all $y \in Y^{\times}$ and $n>M$
$$(\log_p \langle y \rangle)^n   \frac{(1-s)^n}{n!} \frac{1}{p^N} \approx 0 \text{ and }
\binom{1-s}{n} (\langle y \rangle -1)^{n} \frac{1}{p^N} \approx 0$$ This implies that a hyperfinite summation over $n=0,\dots, M$ suffices in \eqref{kubota-logp} and \eqref{kubota-exp}. $M$ depends only on $N$ and not on $s$. The above formulas \eqref{kubota}, \eqref{kubota-logp} and \eqref{kubota-exp} then follow from Lemma \ref{convergence} and do not depend on the choices of $Y$, $N$ and $M$.\\

The nonstandard representation with hyperfinite sums permits convenient computations. For example, 
the residue of $\zeta_{p,1}(s)$ at $s=1$ can be easily computed: $\int_{Y^{\times}} d\nu = \frac{(p-1)p^{N-1}}{p^N} = 1- \frac{1}{p}$. All other branches, i.e. $i \neq 1 \mod (p-1)$, are holomorphic since the integral vanishes at $s=1$:
 $\int_{Y^{\times}} \omega(y)^{1-i} d\nu = \frac{p^{N-1}}{p^N} (\omega(1)^{1-i}+\omega(2)^{1-i}+\dots+\omega(p-1)^{1-i})=0$. Another example is contained in section \ref{mascheroni}.\\
 

The above nonstandard formulas can be reconverted into standard expressions. For example, \eqref{kubota} gives

\[ \zeta_{p,i}(s) =  \frac{(-1)^{i-1}}{s-1} \lim_{n \rightarrow \infty} \sum_{\substack{m=1 \\ p\, \nmid\, m}}^{p^n}  \omega(m)^{1-i} \langle m \rangle^{1-s} \frac{1}{p^n}\]

\subsection{Measures and $p$-adic L-functions}
\label{measureszeta}
Let $p\neq 2$ and $K=\cc_p$. It is well known that regularized Bernoulli {\em measures} on $\zp^{\times}$ can be used to construct
$p$-adic L-functions.
First, the Haar and also the Bernoulli distributions \cite{koblitz} are unbounded in a standard sense.
The Bernoulli distribution $\mu_1$ on $X=\zp$ is defined by $\mu_1(a+p^n \zz_p) = \frac{a}{p^n} - \frac{1}{2}$ for $0\leq a < p^n$. It can be regularized by the following transformation: $\mu(U)=\mu_1(U)-2 \mu_1(\frac{1}{2} U)$ for clopen subsets $U \subset \zp$. The measure $\mu$ coincides with $\mu_{1,\frac{1}{2}}$ in \cite{koblitz} II.5 and $E_2$ in \cite{washington} 12.2. We compute the values of $\mu$. Let $0\leq a < p^n$. We have $\mu(a+p^n \zz_p) = \frac{a}{p^n} - \frac{1}{2} - 2( \{ \frac{a/2}{p^n} \} - \frac{1}{2})$ where $\{\ \}$ denotes the fractional part of a $p$-adic number. We conclude that $\mu(a+p^n \zz_p) = \frac{1}{2}$ if $a$ is even. If $a=2b+1$ is an odd integer, then $\mu(a+ p^n \zz_p) = \frac{1}{p^n} -2 \{ \frac{1/2}{p^n} \} + \frac{1}{2}$. Now the $p$-adic expansion $\frac{1}{2} = \frac{p+1}{2} + \frac{p-1}{2} p + \frac{p-1}{2} p^2 + \dots$ shows that
 $2 \{ \frac{1/2}{p^n} \} = (p+1)p^{-n} + (p-1)p^{-n+1} + \dots + (p-1)p^{-1} = p^{-n}+1$ and hence $\mu(a+p^n \zp)= - \frac{1}{2}$ in that case. This shows that $\mu$ only takes the values $\pm \frac{1}{2}$ on clopen balls. \\
 
The $p$-adic zeta functions $\zeta_{p,i}(s)$, and more general $p$-adic L-functions $L_p(s,\chi)$ for arbitrary Dirichlet characters $\chi$, can be defined as a $p$-adic Mellin transform of a measure on $X^{\times}=\zp^{\times}$ (see for example \cite{washington} ch.\ 12). The above regularized Bernoulli measure $\mu$ can be used to define $p$-adic L-functions $L_p(s,\omega^k)$ for powers of the Teichmüller characters $\omega$ and one has

\[ \zeta_{p,i}(s) = L(s,\omega^{1-i}) = \frac{-1}{1-\omega(2)^{1-i} \langle 2 \rangle^{1-s}} \int_{X^{\times}} \omega(x)^{-i} \langle x \rangle^{-s} d\mu \]

Now let $N \in \nsn$ be infinite, $Y={^*\zp}/p^N {^*\zp}$ and $Y^{\times} = ({^*\zp}/p^N {^*\zp})^{\times}$ hyper\-finite spaces with a measure defined by $\nu(\{a+p^N {^*\zp}\})=\frac{(-1)^a}{2}$ for $0\leq a < p^N$ (see Proposition \ref{hyperfinite} and subsequent examples).  Theorem \ref{hyperfinite-rep} shows that the measure space $(X,B(X),\mu)$ is the push-down of $(Y,{^*\mathcal{P}}(Y),\nu)$ via the standard-part map. Although there are obviously measurable subsets of $X$ and $Y$ with zero measure (e.g.\ the union of two balls or points with measure $\frac{1}{2}$ and $-\frac{1}{2}$), the weight functions $N_{\mu}$ and $N_{\nu}$ have the constant value $|\frac{1}{2}|=1$. An internal function $F: Y \rightarrow {^*\cc_p}$ is $S$-integrable if and only if $\|F\|_{\nu}$ is finite, i.e. if $|F(y)|$ is bounded by a standard real number.  \\

As above, the continuous function $f(x)=\langle x \rangle^{s}$ on $X^{\times}$ can be lifted to a function $F(y)=\langle y \rangle^{s}$ on
$Y^{\times}$
such that $f(\st_Y (y))=\st_K( F(y))$ for $y \in Y^{\times}$.  
One may can take the $^*$-polynomial $F(y)= \sum_{n=0}^M \binom{1-s}{n} (\langle y \rangle -1)^{n}$ where $M \in \nsn$ depends on $N$.
We obtain a nonstandard representation of $\zeta_{p,i}(s)$:

\begin{align} \zeta_{p,i}(s) = & \  \frac{-1}{1-\omega(2)^{1-i} \langle 2 \rangle^{1-s}} \st_K \left( \int_{Y^{\times}} \omega(y)^{-i} \langle y \rangle^{-s} d\nu \right) \nonumber \\
= & \ \frac{-1}{1-\omega(2)^{1-i} \langle 2 \rangle^{1-s}} \st_K \left( \sum_{\substack{m=0 \\ p\, \nmid\, m}}^{p^N} \omega(m)^{-i} \langle m \rangle^{-s} \frac{(-1)^m}{2} \right) 
\label{zeta}
\end{align}

A standard interpretation of this integral is the following series:

\[ \zeta_{p,i}(s) = \frac{1}{2(1-\omega(2)^{1-i} \langle 2 \rangle^{1-s})} \sum_{n=1}^{\infty} \left( \sum_{\substack{m=p^{n-1} \\ p\, \nmid\, m,\ m \text{ odd} }}^{p^n} \omega(m)^{-i} \langle m \rangle^{-s} -
\sum_{\substack{m=p^{n-1} \\ p\, \nmid\, m,\ m \text{ even} }}^{p^n} \omega(m)^{-i} \langle m \rangle^{-s} \right)
\]

For $s=k$ and $i=k \mod p-1$ one has

\[ \zeta_{p,i}(k) = L_p(k,\omega^{1-k}) = \frac{-1}{1-2^{1-k}}  \cdot \lim_{n \rightarrow \infty} \sum_{m=1}^{p^n} \frac{(-1)^m}{2} m^{-k}  \]

Such Dirichlet series expansions were proved by D. Delbourgo in \cite{delbourgo2006dirichlet} using $p$-adic fractional derivations. 
$p$-adic Euler products and series expansions for arbitrary Dirichlet characters are treated in his work \cite{delbourgo2009} and it might be possible to obtain similar results with nonstandard methods, but we leave this for future work.

\subsection{$p$-adic Euler-Mascheroni constant}
\label{mascheroni}

Let $p \neq 2$, $K=\cc_p$, $\st=\st_K$ the standard-part map and $\zeta_{p,\mathbf{1}}(s)= \frac{1-\frac{1}{p}}{s-1} + \gamma_p + \dots$ the Laurent expansion of the $\mathbf{1}$-branch of the $p$-adic zeta function around $s=1$. The coefficient $\gamma_p$ is called the $p$-adic {\em Euler-Mascheroni} constant. We derive different formulas for that constant.\\

\noindent A) First, $\gamma_p$ can be computed with the Kubota-Leopoldt zeta function and the internal Haar measure $\nu$ (see section \ref{klzeta}). By \eqref{kubota-org} above, we find that $\gamma_p$ is the derivative of 
  $-\int_{X^{\times}} \langle x \rangle^{1-s} dx$
 at $s=1$.  We use \eqref{kubota}, set $1-s=p^N \approx 0$ and obtain the following nonstandard formula for $\gamma_p$ :
 
 \begin{equation} \gamma_p = \st \left( \frac{1}{-p^N} \left (  \int_{Y^{\times}} \langle y \rangle^{p^N} d\nu \ - (1-\frac{1}{p}) \right) \right) 
 = \st \left( \frac{-1}{p^N} \left (  \frac{1}{p^N} \sum_{\substack{m=1 \\ p\, \nmid m}}^{p^N} \langle m \rangle^{p^N}  - (1-\frac{1}{p}) \right) \right) 
 \label{first}
 \end{equation}
We have $\langle m \rangle^{p^N} = \omega^{-1}(m) m^{p^N}$ and the sums of $p^N$-th powers can be computed with generalized Bernoulli numbers:
\[ \sum_{\substack{m=1 \\ p\, \nmid m}}^{p^N} \langle m \rangle^{p^N} = \sum_{m=1}^{p^N} \omega^{-1}(m) m^{p^N} = p^N B_{p^N,\omega^{-1}} + \frac{p^N}{2} (p^N)^2 B_{p^N-1,\omega^{-1}} + \dots \]

By the Theorem of von Staudt-Clausen  for generalized Bernoulli numbers \cite{leopoldt}, $|B_{k,\omega^{-1}}|_p \leq p$. Therefore
\[ \st \left( \frac{1}{p^N} \sum_{m=1}^{p^N} \omega^{-1}(m) m^{p^N}  \right) = B_{p^N,\omega^{-1}} \]
Hence it follows that
\begin{equation} \gamma_p = \st \left( \frac{-1}{p^N} \left( B_{p^N,\omega^{-1}} - (1-\frac{1}{p}) \right) \right) 
= \lim_{n \rightarrow \infty} \frac{1}{p^n} \left( 1-\frac{1}{p} - B_{p^n,\omega^{-1}}   \right)
\label{knospe}
\end{equation}
We remark that a similar formula can be found in \cite{koblitz1978}. \\

\noindent B) Alternatively, one can use the expansion \eqref{kubota-logp} which gives

\begin{equation} \gamma_p = \st \left(- \int_{Y^{\times}} \log_p \langle y \rangle  d\nu \right) = \st \left( - \frac{1}{p^N}  \sum_{\substack{m=1 \\ p\, \nmid m}}^{p^N} \log_p(m)  \right)
\label{kubota-log}
\end{equation}

One has $\displaystyle\sum_{\substack{m=1 \\ p\, \nmid m}}^{p^N} \log_p(m)=  \log_p \left( \displaystyle\prod_{\substack{m=1 \\ p\, \nmid m}}^{p^N} m \right) = \log_p (\Gamma_p(p^N))$, where $\Gamma_p$ is the $p$-adic gamma function. Since $\Gamma_p(0)=1$ and $\Gamma_p(x+1)=-\Gamma_p(x)$ for $|x|<1$, we get $\st (\frac{1}{p^N} \log_p (\Gamma_p(p^N)) )= (\log_p \Gamma_p)'(0) = \frac{\Gamma_p '(0)}{\Gamma_p(0)} = -\Gamma_p'(1)$. $\Gamma_p$ interpolates the factorial with the factors dividing $p$ removed and this yields the following formula:
\begin{equation} \gamma_p = \Gamma_p'(1) = \st \left( \frac{1}{p^N} \left(\Gamma_p(p^N+1) -(-1) \right) \right) = \st \left( \frac{1}{p^N}  \left( \frac{p^N !}{p^{N-1}! \ p^{p^{N-1}}} +1 \right) \right)
\label{gamma-fak}
\end{equation}
The standard version of \eqref{gamma-fak} can be found in \cite{schikhof} 36.A. \\

\noindent C) Next, we use the expansion \eqref{kubota-exp} and set $1-s=p^N$. This gives

\begin{align} \gamma_p = & - \st \left( \frac{1}{p^N} \int_{Y^{\times}} \sum_{n=1}^M \binom{p^N}{n} (\langle y \rangle -1)^{n} \ d\nu \right) \nonumber \\ 
= &  - \st \left( \int_{Y^{\times}} \sum_{n=0}^M \frac{(-1)^{n+1}}{n} \sum_{j=0}^n \binom{n}{j} \left(\frac{y} {\omega(y)}\right)^{j} (-1)^{n-j} \ d\nu \right) \nonumber \\
= &  - \st \left( \sum_{n=0}^M \frac{(-1)^{n+1}}{n} \sum_{j=0}^n \binom{n}{j}  (-1)^{n-j} \int_{Y^{\times}}  \omega^{-j}(y) y^{j}  \ d\nu \right) \nonumber \\
= &  - \st \left( \sum_{n=0}^M \frac{1}{n} \sum_{j=0}^n \binom{n}{j}  (-1)^{j+1} B_{j,\omega^{-j}} \right)
\label{newknospe}
\end{align}

We also state a standard version of the series \eqref{newknospe} which converges very fast.
\[ \gamma_p = \lim_{m \rightarrow \infty}  \left( \sum_{n=0}^m \frac{1}{n} \sum_{j=0}^n \binom{n}{j}  (-1)^{j+1} B_{j,\omega^{-j}} \right) \]

\noindent D) One may also use the regularized Bernoulli measure and its nonstandard representation $\nu$ as in section \ref{measureszeta}. By \eqref{zeta} above,

\[ \zeta_{p,\mathbf{1}}(s) =  \frac{-1}{1- \langle 2 \rangle^{1-s}} \st \left( \int_{Y^{\times}} \frac{1}{y} \langle y \rangle^{1-s} d\nu \right) \]

We compute $\gamma_p$ by setting $1-s=p^N$ and get

\begin{align*} \gamma_p = & \ \st \ \frac{-1}{1- \langle 2 \rangle^{p^N}}    \left( \int_{Y^{\times}} \frac{1}{y} \langle y \rangle^{p^N} d\nu \right) - \frac{1-\frac{1}{p}}{-p^N}  \\
 = & \ \st \  \frac{1}{1-\langle 2 \rangle^{p^N}}  \left(  \sum_{\substack{m=1 \\ p\, \nmid m}}^{p^N} \frac{(-1)^{m+1}}{2m} \langle m \rangle^{p^N} + (1-\frac{1}{p})(1-\langle 2 \rangle^{p^N}) \frac{1}{p^N} \right)
 \end{align*}

Since $\langle 2 \rangle^{p^N} = \exp_p(p^N \log_p\langle 2 \rangle) = 1 + p^N \log_p\langle 2 \rangle + \frac{1}{2} p^{2N} \log_p\langle 2 \rangle^2 + \dots$, we conclude
\begin{equation} \gamma_p = \ \st  \ \frac{1}{1-\langle 2 \rangle^{p^N}}  \left(  \sum_{\substack{m=1 \\ p\, \nmid m}}^{p^N} \frac{(-1)^{m+1}}{2m} \langle m \rangle^{p^N} - (1-\frac{1}{p})\log_p\langle 2 \rangle (1 + \frac{1}{2} p^N  \log_p\langle 2 \rangle  )  \right) 
\label{delb}
\end{equation}

A similar formula for $\gamma_p$ modulo $p^n \zp$ was proved by D. Delbourgo \cite{delbourgo2009} 2.7. \\

\noindent E) Yet another formula for $\gamma_p$ follows from the construction of $p$-adic L-functions in \cite{washington} 5.11. One has the following expansion around $s=1$:
\begin{align} (s-1)\cdot L_p(s,\mathbf{1})  
= &  \frac{1}{p}  \sum_{a=1}^{p-1} \exp_p( (1-s) \log_p \langle a \rangle) \sum_{j=0}^{\infty} \binom{1-s}{j} B_j \frac{p^j}{a^j} \nonumber
\\ = & \frac{p-1}{p} + \frac{1}{p} \left ( \sum_{a=1}^{p-1} - \log_p \langle a \rangle - B_1 \frac{p}{a} + \frac{1}{2} B_2 \frac{p^2}{a^2} - \frac{1}{3} B_2 \frac{p^3}{a^3} \pm \dots \right) (s-1) + \dots 
\label{wash}
\end{align}
The linear term of the expansion is equal to $\gamma_p$ and this series converges very fast. \\

We did a number of calculations using the mathematics software {Sage} \cite{sage} to provide additional numerical evidence for \eqref{first}, \eqref{knospe}, \eqref{kubota-log}, \eqref{newknospe}, \eqref{delb}, \eqref{wash}.
  
For example, the formulas give the following values (compare \cite{delbourgo2009}):
\begin{align*} \gamma_3 =\ & 2\cdot3 + 2\cdot 3^2 + 3^3 + 2\cdot 3^4 + 3^5 + 2\cdot 3^6 + 2\cdot 3^7 + 2\cdot 3^8 + O(3^{10}) \\
 \gamma_5 =\ & 5 + 3\cdot5^3 + 2\cdot5^5 + 3\cdot5^6 + 4\cdot5^7 + 5^8 + 2\cdot5^9 + O(5^{10})\\
 \gamma_7 =\ & 5 + 2 \cdot 7 + 4\cdot 7^2 + 6\cdot 7^3 + 2\cdot7^4 + 6\cdot7^6 + 2\cdot 7^7 + 7^9 + O(7^{10}) \\
 \gamma_{11} =\ & 1 + 10 \cdot 11 + 2 \cdot 11^2 + 11^3 + 5 \cdot 11^4 + 5 \cdot 11^5 
 + 4 \cdot 11^6 
 + 5 \cdot 11^7 + O(11^8) \\
\gamma_{13} =\ & 4\cdot 13 + 7\cdot13^3 + 8\cdot13^4 + 7\cdot13^5 + 6\cdot13^6 + 4\cdot13^7+9\cdot13^8 +O(13^9)
\end{align*}

\section*{Acknowledgments}
The author wishes to express his thanks to Daniel Delbourgo for drawing the author's attention to 
Dirichlet series expansions of $p$-adic L-functions and to the $p$-adic Euler-Mascheroni constant. Furthermore, he wants to thank Jörn Steuding for pointing him to the $p$-adic Gamma function.

\bibliography{nsabib}

\bibliographystyle{plain}

\end{document}